\documentclass[a4paper,10pt]{article}

\usepackage{graphicx}
\usepackage{amsmath}
\usepackage{theorem}
\usepackage{amscd}
\usepackage{amssymb}
\usepackage[T1]{fontenc}
\usepackage{epsfig}
\usepackage{psfrag}

\newenvironment{demo}{\noindent \textbf{Proof :\\}}{q.e.d.}
\newenvironment{preuve}{\noindent \textbf{Proof :\\}}{q.e.d.}

\newtheorem{thm}{Theorem}[section]
\newtheorem{cor}[thm]{Corollary}
\newtheorem{lem}[thm]{Lemma}
\newtheorem{prop}[thm]{Proposition}
\newtheorem{defn}[thm]{Definition}
\newtheorem{rem}[thm]{Remark}
\newtheorem{ex}[thm]{Example}
\theoremstyle{break}

\newcommand{\norm}[1]{\left\Vert#1\right\Vert}
\newcommand{\abs}[1]{\left\vert#1\right\vert}

\newcommand{\set}[1]{\left\{#1\right\}}
\newcommand{\Forall}[2]{\forall \, #1 \in #2, \:}
\newcommand{\Exists}[2]{\exists \, #1 \in #2, \:}

\newcommand{\restric}[1]{\vert_{#1}}

\newcommand{\Real}{\mathbb R}
\newcommand{\Complex}{\mathbb C}
\newcommand{\Integer}{\mathbb Z}
\newcommand{\Natural}{\mathbb N}
\newcommand{\Sphere}{\mathbb S}

\newcommand{\To}{\longrightarrow}
\newcommand{\U}{\mathcal{U}}
\newcommand{\E}{\mathcal{E}}
\newcommand{\V}{\mathcal{V}}
\newcommand{\G}{\mathcal{G}}
\newcommand{\A}{\mathcal{A}}
\newcommand{\B}{\mathcal{B}}
\newcommand{\C}{\mathcal{C}}
\newcommand{\M}{\mathcal{M}}
\newcommand{\X}{\mathcal{X}}
\newcommand{\Y}{\mathcal{Y}}

\DeclareMathOperator{\length}{length}

\DeclareMathOperator{\Ric}{Ric}

\DeclareMathOperator{\Card}{Card}
\DeclareMathOperator{\vol}{vol}
\DeclareMathOperator{\supp}{supp}
\DeclareMathOperator{\Id}{Id}
\DeclareMathOperator{\grad}{grad}

\begin{document}


\title{\bfseries{Weighted Sobolev inequalities and Ricci flat manifolds.}}
\author{Vincent Minerbe}
\date{\today}

\maketitle

\abstract{In this paper, we prove a weighted Sobolev inequality and a Hardy inequality on manifolds with nonnegative Ricci curvature 
satisfying an inverse doubling volume condition. It enables us to obtain rigidity results for Ricci flat manifolds.}


\section*{Introduction}

Since the eighties and particularly \cite{BKN}, it is well known that Ricci flat manifolds with maximal volume growth enjoy nice rigidity 
properties. Indeed, if $M^n$, $n\geq 4$, is such a manifold, with curvature tensor $R$, there exists a constant $\epsilon$ such that $M$ is 
flat as soon as $\int_M \abs{R}^{n/2} dvol < \epsilon$ ; $\epsilon$ depends on $n$ and on a lower bound on the volume growth. Furthermore, 
in case the curvature only satisfies $\int_M \abs{R}^{n/2} dvol < \infty$, it has faster-than-quadratic decay, that is 
$R = O(r_o^{-2-\delta})$, where $r_o$ is the geodesic distance to any point $o$ in $M$ ; here, $\delta$ is an explicit positive constant. 
These facts stem from a Sobolev inequality which is no longer true if the volume growth is not maximal. Now what happens in this case ? One 
result in this direction is the following theorem, by Jeff Cheeger and Gang Tian \cite{CT} : a four-dimensional complete Ricci flat manifold 
with curvature in $L^2$ has quadratic curvature decay. Their proof is based on the Gauss-Bonnet-Chern formula and Cheeger-Gromov theory.

Our aim is to present a different approach, relying on weighted Sobolev and Hardy inequalities. Unlike J. Cheeger and G. Tian, we still make 
an assumption on the volume growth, and this enables us to generalize the rigidity results which were previously known. Given a point $o$, 
we will consider weights involving the function $\rho_o :\, t \mapsto \frac{t^n}{V(o,t)}$, where $V(o,t)$ is the volume of the geodesic 
ball $B(o,t)$ centered in $o$ and of radius $t$. Our work leads to the following.

\begin{thm}[Flatness criterion]
Let $M^n$, $n \geq 4$, be a connected complete Ricci-flat manifold. Assume there
exists $o$ in $M$, $\nu > 1$ and  $C_o >0$ such that 
$$
\forall t_2 \geq t_1 > 0, \  \frac{V(o,t_2)}{V(o,t_1)} \geq C_o \left( \frac{t_2}{t_1} \right)^\nu.
$$
Then there is a constant $\epsilon_1(n,C_o,\nu)$ such that $M$ is flat as soon as
$$
\sup_M (\abs{R} r_o^2) < \epsilon_1(n,C_o,\nu).
$$
If $\nu>2$, there is also a constant $\epsilon_2(n,C_o,\nu)$ such that $M$ is flat as soon as 
$$
\int_M \abs{R}^{\frac{n}{2}} \rho_o(r_o) dvol < \epsilon_2(n,C_o,\nu).
$$
As a result, in both cases, $M$ is the normal bundle of a compact totally geodesic submanifold, which is (finitely) covered by a flat torus.   
\end{thm}

\begin{thm}[Curvature decay]
Let $M^n$, $n \geq 4$, be a connected complete Ricci-flat manifold. Assume there
exists $o$ in $M$, $\nu > 2$ and  $C_o >0$ such that 
$$
\forall t_2 \geq t_1 > 0, \  \frac{V(o,t_2)}{V(o,t_1)} \geq C_o \left( \frac{t_2}{t_1} \right)^\nu.
$$
and 
$$
\int_M \abs{R}^{\frac{n}{2}} \rho_o(r_o) dvol < + \infty.
$$
Then $M$ has quadratic curvature decay. Furthermore, if $\nu > 4\frac{n-2}{n-1}$, $M$ has faster-than-quadratic curvature decay and 
thus has finite topological type. 
\end{thm}

The assumption 
\begin{equation}\label{volboules1}
\forall t_2 \geq t_1 > 0, \  \frac{V(o,t_2)}{V(o,t_1)} \geq C_o \left( \frac{t_2}{t_1} \right)^\nu
\end{equation}
implies 
\begin{equation}\label{volboules1'}
\forall R\geq 1, \  V(o,R) \geq C_o V(o,1) R^\nu
\end{equation}
and follows from    
\begin{equation}\label{volboules2}
\exists A_o, B_o >0, \, \forall R\geq 1, \  A_o R^\nu \leq V(o,R) \leq B_o R^\nu.
\end{equation}
Note that Bishop theorem ensures $\nu \leq n$. This hypothesis yields the analytical tools we need. Indeed, we prove that on a complete 
connected manifold $M^n$, $n\geq3$, with nonnegative Ricci curvature and satisfying (\ref{volboules1}), the following weighted Sobolev 
inequality holds :
\begin{equation}\label{introsobo}
\Forall{f}{C^\infty_c(M)} \left(\int_M  \abs{f}^{\frac{2 n}{n-2}} \rho_o(r_o)^{-\frac{2}{n-2}} dvol \right)^{\frac{n-2}{n}} 
\leq S \int_M \abs{df}^2 dvol. 
\end{equation}
In other terms, the completion $H^1_0 (M)$ of $C^\infty_c(M)$ for the norm $\norm{d.}_{L^2(M,vol)}$ can be 
continuously injected into $L^{\frac{2 n}{n-2}}\left(M,\rho_o(r_o)^{-\frac{2}{n-2}} vol \right)$. Such a manifold 
also satisfies the Hardy inequality 
\begin{equation}\label{introhardy}
\Forall{f}{C^\infty_c(M)} \int_M  \abs{f} r_o^{-1}  dvol \leq H \int_M \abs{df} dvol.
\end{equation}
The constants $S$ and $H$ we find depend only on $n$, $\nu$ and $C_o$.
Now, the curvature of a Ricci flat manifold obeys a nonlinear elliptic equation. When used appropriately, the inequalities (\ref{introsobo}) and (\ref{introhardy}) yield estimates on the solutions of this equation, and our theorems follow. In this article, we will give a few other applications of the weighted Sobolev inequality.

\vspace{\baselineskip}

The paper is organized as follows.

In the first section, we develop a discretization technique aimed at patching local Sobolev inequalities together. It is based upon ideas and 
methods of A. Grigor'yan and L. Saloff-Coste \cite{GSC}. Given a convenient covering of a manifold, if we assume some discrete inequality on 
a graph which is naturally associated to the covering, we are able to deduce a global Sobolev inequality from a local one (theorem 
\ref{recollsobodir}).

In the second section, we explain how to apply this abstract technique in the setting of manifolds with nonnegative Ricci curvature and satisfying (\ref{volboules1}), so as to
obtain a weighted Sobolev inequality and a Hardy inequality. Note we could replace the nonnegativity of the Ricci curvature by two of its
consequences : the doubling volume condition and the scaled Poincaré inequality on balls. In \cite{Gril}, G. Grillo proves weighted inequalities in the context of homogeneous spaces and indeed, in the case $\nu=n$, the Hardy inequality follows from this work : nevertheless, it
should be stressed that our approach is basically different and, in particular, does not require a uniform estimate on the volume of balls ; 
apart from the doubling volume condition and
the scaled Poincaré inequality (which are classical assumptions for such problems), the only measure theoretic assumption we need is the 
estimate (\ref{volboules1}) which is some kind of inverse doubling volume condition \emph{around  one point}. An important step in our proof 
could be singled out : the following result gives a sufficient condition for a manifold to satisfies the so called RCA property 
(Relatively Connected Annuli) and should be compared with proposition 4.5 of \cite{HK} (which, in our context, would require the volume 
growth to satisfy a uniform euclidian estimate from below). 
\begin{prop}[RCA]
Let $M$ be a connected complete riemannian manifold, satisfying the doubling volume property
$$
\Forall{x}{M} \forall \, R>0, \,  V(x,2R) \leq C_D V(x,R),
$$
the scaled $L^p$ Poincaré inequality centered in some point $o$ in $M$ 
$$
\Forall{f}{C^\infty_c(M)} \forall \, R>0, \,  \int_{B(o,R)} \abs{f - f_{B(o,R)}}^p dvol \leq C_P R^p \int_{B(o,R)} \abs{df}^p dvol
$$
and the inverse doubling volume condition centered in $o$
$$
\forall R_2 \geq R_1 >0, \  \frac{V(o,R_2)}{V(o,R_1)} \geq C_o \left( \frac{R_2}{R_1} \right)^\nu
$$
with $\nu > p$. Here, $C_D \geq 1$, $p \geq 1$, $C_P >0$, $C_o >0$. Then there exists $\kappa_0 >0$ such that for $R > 0$, 
if $x,y$ are two points in $S(o,R)$, there is a path from $x$ to $y$ which remains inside $B(o,R)\backslash B(o,\kappa_0^{-1} R)$. Moreover,
it is possible to find an explicit constant, in terms of $p,C_D,C_P,C_o,\nu$.  
\end{prop}
Let us say a few words about this proposition. Cheeger-Gromoll theorem implies that in our setting, $M$ has only one end. A result from \cite{LT} 
(with \cite{And}) implies that, for large $R$, the intersection of the only unbounded component of $M\backslash B(o,R)$ with any annulus 
$A(R,R+r)$, $r>0$, is connected. But it says nothing about the behaviour of the bounded components of $M\backslash B(o,R)$. What we proved is
that, in a sense, these bounded components have at most linear growth. Moreover, we give an explicit estimate of this growth, which is
important for our purpose. 

In the third section, we investigate the properties of Schrödinger operators $\overline{\Delta} + V$ that can be deduced from our weighted 
Sobolev inequality. Here, $\overline{\Delta}$ is the Bochner laplacian on some euclidian vector bundle and $V$ is a field of symmetric 
endomorphisms. In particular, we prove that integral assumptions on the potential ensure the kernel is trivial (theorem \ref{annulation}). 
We obtain various technical estimates and also introduce a good space of sections $\psi$ such that the equation 
$
(\overline{\Delta} + V)\sigma = \psi
$    
has a bounded solution $\sigma$ (\ref{inv3}). This section can be seen as a toolbox.

In the fourth and last section, we point out some applications. Let us denote by $S_{o}(M)$ (resp. $H_o(M)$) the best constant $S$ (resp. $H$) in (\ref{introsobo}) (\ref{introhardy}). We define the "Sobolev-curvature" invariant 
$$
\mathcal{SC} (M^n):= \inf_{o \in M} \left[ S_{o}(M) \left(\int_M \abs{R}^{\frac{n}{2}} \frac{r_o^n}{V(o,r_o)} dvol \right)^{\frac{2}{n}} \right]
$$
and the "Hardy-curvature" invariant
$$
\mathcal{HC} (M^n):= \inf_{o \in M} \left[ H_{o}(M) \sup_M (\abs{R} r_o^2) \right],
$$
with the convention $0.\infty=\infty$. First, we generalize the work of G. Carron \cite{Car1} about
$L^2$-cohomology and obtain in particular the
\begin{thm}[$L^2$-cohomology]
Let $M^n$, $n \geq 3$, be a connected complete riemannian manifold such that $\mathcal{SC} (M)$ is finite. 
Then the $L^2$ Betti numbers of $M$ are finite. Moreover, $\mathcal{H}_{L^2}^1(M) = \set{0}$ and, for $k \geq 2$, there exists a positive 
universal constant $\epsilon(n,k)$ such that if
$
\mathcal{SC}(M) < \epsilon(n,k),
$
then $\mathcal{H}_{L^2}^k(M) = \set{0}$. 
\end{thm}
In case $M$ has nonnegative Ricci curvature and satisfies (\ref{volboules2}), this means $L^2$ Betti numbers are finite as soon as $\int_M
\abs{R}^{\frac{n}{2}} r_o^{n-\nu} dvol < \infty$ ; in \cite{Car2}, G. Carron required $\int_M \abs{R}^{\frac{\nu}{2}} dvol < \infty$. 
These are close assumptions, but, for instance, if the curvature decays quadratically (that is $\mathcal{HC} (M) < \infty$, in this setting), 
ours is weaker.  Our work also provides explicit bounds on the $L^2$ Betti numbers. Then we study Ricci flat manifolds and prove the 
following rigidity theorems, which imply the results announced above.  
\begin{thm}[Flatness criterion]
If $M^n$, $n \geq 4$,  is a connected complete Ricci-flat manifold, there are  
universal positive constants $\epsilon_1(n)$ and $\epsilon_2(n)$ such that if 
$\mathcal{SC} (M) < \epsilon_1(n)$ or $\mathcal{HC} (M) < \epsilon_2 (n)$,
then $M$ is flat. 
\end{thm}
\begin{thm}[Curvature decay]
Let $M^n$, $n \geq 4$, be a connected complete Ricci-flat manifold. If $\mathcal{SC}(M)$ is finite, 
then $M$ has quadratic curvature decay. If moreover there exists $\nu > 4\frac{n-2}{n-1}$ and $A_o > 0$ such that $V(o,R) \geq A_o R^\nu$ 
for large $R$, then the curvature decays like $r_o^{- \frac{(\nu-2)(n-1)}{n-3}}$ ; in particular, $M$ has finite topological type.
\end{thm}
Finally, we give some examples where this rate of decay is the correct one.

\vspace{\baselineskip}

$\bf{Acknowledgements}$. I would like to thank Gilles Carron for his remarks, suggestions, questions, and for his patience.

\newpage

\tableofcontents

\newpage

\section{Discretization and Sobolev inequalities.}

\subsection{How to patch local Sobolev inequalities together.}

The aim of this paragraph is to explain how to patch local Sobolev inequalities so as to obtain a global one. In \cite{GSC}, 
A. Grigor'yan and L. Saloff-Coste introduce a discretization procedure enabling them to handle Poincaré inequalities.
We generalize their ideas in two ways : we use integral inequalities for different measures and we consider general Sobolev-type 
inequalities. 

Here, $M$ is a smooth riemannian manifold (Lipschitz would be sufficient), and we introduce two Borel measures $\lambda$, $\mu$ 
on it. For us, it will be crucial to cope with both of them at the same time. Let us introduce the necessary vocabulary.

\begin{defn}\label{bonrec}
Let $A \subset A^\sharp$ be two subsets of $M$. A family $\U=(U_i,U_i^*,U_i^\sharp)_{i \in I}$ 
consisting of subsets of $M$ having finite measure with respect to $\lambda$ and $\mu$ is said to be a good covering of $A$ in $A^\sharp$ 
if the following is true.
\begin{itemize}
\item[(i)] There is a Borel subset $E$ of $A$ with $\lambda(E)=\mu(E)=0$, such that 
$A \backslash E \subset \bigcup_{i} U_i \subset \bigcup_{i} U_i^\sharp \subset A^\sharp$;\
\item[(ii)] $\Forall{i}{I} U_i \subset U_i^*\subset U_i^\sharp$ ;\
\item[(iii)] There exists a constant $Q_1$ such that for each $i_0 \in I$,
 $$\Card \set{i \in I / U_{i_0}^\sharp \cap U_{i}^\sharp \not= \emptyset} \leq Q_1;$$\
\item[(iv)] For every $(i,j) \in I^2$ satisfying $\overline{U_{i}} \cap \overline{U_{j}} \neq \emptyset$, there is an element
$k(i,j)$ such that $$U_{i} \cup U_{j} \subset U_{k(i,j)}^*;$$\
\item[(v)] There exists a constant $Q_2$ such that for every $(i,j) \in I^2$, if $\overline{U_{i}} \cap \overline{U_{j}} \neq \emptyset$, then
 $$\lambda(U_{k(i,j)}^*) \leq Q_2 \min{(\lambda(U_{i}),\lambda(U_{j}))}$$ and $$\mu(U_{k(i,j)}^*) \leq Q_2 \min{(\mu(U_{i}),\mu(U_{j}))}.$$
\end{itemize}
\end{defn}

Given a Borel set $U$ with finite and nonzero $\lambda$-measure \ and a $\lambda$-integrable function $f$, we will denote by $f_{U,\lambda}$ the mean value of $f$ on $U$ with respect to the measure $\lambda$ :
$$
f_{U,\lambda} = \frac{1}{\lambda(U)} \int_U f d\lambda.
$$

One can associate to every good covering $\U$ a weighted graph $(\G,m_\lambda)$ : its set of vertices is  $$\V = I$$ and its set of edges
is $$\E = \set{\set{i,j} \subset \V / i \not= j, \, \overline{U_{i}} \cap \overline{U_{j}} \neq \emptyset} ;$$ $\V$ and $\E$ are given 
measures, both of which will be denoted by $m_\lambda$, and they are defined by 
$$\Forall{i}{\V} m_\lambda(i) = \lambda (U_i)$$
and
$$\Forall{(i,j)}{\E} m_\lambda(i,j)=\max(m_\lambda(i),m_\lambda(j)).$$ 

\begin{rem}\label{graphepondere}
In what we call a graph, there is at most one edge between two given vertices. So, if there is an edge between two vertices $i$ and $j$, we
will call it $(i,j)$. For us, a weighted graph will always consist of a
$\sigma$-finite measure on the set of vertices $\V$ and of a $\sigma$-finite measure on the set of edges $\E$, which we give the same name $m$ and which are related by
$$
\Forall{(i,j)}{\E} m(i,j)=\max(m(i),m(j)).
$$ 
\end{rem}

Now, we introduce three kinds of inequalities : the discrete estimates (the second and third) will enable us to patch the continuous ones
(the first) together. 

\begin{defn}\label{sc}
Suppose $k \in ]1,\infty]$ and $1\leq p < k$. We will say that a good covering $\U$ satisfies a continuous $L^p$ Sobolev inequality  of order 
$k$ with respect to the pair of measures $(\lambda,\mu)$ if there exist a constant $S_c$ such that for every $i \in I$, one has 
$$
\Forall{f}{C^\infty(U_i^*)}  \left(\int_{U_i} \abs{f - f_{U_i,\lambda}}^{\frac{p k}{k-p}} d\lambda \right)^{\frac{k-p}{k}} 
\leq S_c \int_{U_i^*} \abs{df}^p d\mu,
$$  
and
$$
\Forall{f}{C^\infty(U_i^\sharp)}  \left(\int_{U_i^*} \abs{f - f_{U_i^*,\lambda}}^{\frac{p k}{k-p}} d\lambda \right)^{\frac{k-p}{k}} 
\leq S_c \int_{U_i^\sharp} \abs{df}^p d\mu.
$$  
\end{defn}

\begin{defn}\label{sdd}
Suppose $k \in ]1,\infty]$ and $1\leq p < k$. We will say that the weighted graph $(\G,m)$ satisfies a discrete $L^p$ Sobolev-Dirichlet
inequality of order $k$ if there exists a constant
$S_d$ such that for every $f \in L^p(\V,m)$, one has 
$$
\left(\sum_{i \in \V} \abs{f(i)}^{\frac{p k}{k-p}} m(i)\right)^{\frac{k-p}{k}} 
\leq S_d \sum_{(i,j) \in \E} \abs{f(i)-f(j)}^p m(i,j).
$$  
\end{defn}

\begin{defn}\label{sdn}
Suppose $k \in ]1,\infty]$ and $1\leq p < k$. We will say that a finite weighted graph $(\G,m)$ satisfies a discrete $L^p$ Sobolev-Neumann 
inequality of order $k$ if there exists a constant $S_d$ such that for every $f \in \Real^\V$, one has 
$$
\left(\sum_{i \in \V} \abs{f(i) - m(f)}^{\frac{p k}{k-p}} m(i)\right)^{\frac{k-p}{k}} 
\leq S_d \sum_{(i,j) \in \E} \abs{f(i)-f(j)}^p m(i,j).
$$  
\end{defn}

\begin{rem}
In this terminology, a $L^p$  Poincaré inequality is nothing but a $L^p$ Sobolev inequality of infinite order.
\end{rem}

\begin{rem}
Of course, one can say that a good covering $\U$ satisfies a discrete Sobolev inequality, by considering the associated weighted graph
$(\G,m_\lambda)$.
\end{rem}

The following theorem is the crucial tool for us.

\begin{thm}\label{recollsobodir}
Suppose $k \in ]1,\infty]$ and $1\leq p < k$. If a good covering $\U$ of $A$ in $A^\sharp$ satisfies  
the continuous $L^p$ Sobolev inequality of order $k$ (\ref{sc}) and the discrete $L^p$ Sobolev-Dirichlet of order $\infty$ (\ref{sdd}), 
then the following Sobolev-Dirichlet inequality is true :
$$\Forall{f}{C^\infty_c(A)} \int_A \left( \abs{f}^{\frac{p k}{k-p}} d\lambda \right)^{\frac{k-p}{k}} 
\leq S \int_{A^\sharp} \abs{df}^p d\mu, 
$$
with 
$$
S= 2^{p-1 + \frac{p}{k}} ((S_c Q_1)^{\frac{k}{k-p}} + S_d Q_2 (2^p S_c Q_1^3)^{\frac{k}{k-p}})^{\frac{k-p}{k}}.
$$
\end{thm}

\begin{rem}
The case where $\lambda = \mu$, $k=\infty$ and $p=2$ was proved by A. Grigor'yan and L. Saloff-Coste in \cite{GSC}. 
\end{rem}

\begin{demo}
We set $q:=\frac{pk}{k-p}$ and consider $f \in C^\infty_c(A)$. 
Thanks to a little convexity, we can write 
\begin{eqnarray*}
\int_A \abs{f}^q d\lambda &\leq& \sum_{i \in \V} \int_{U_i} \abs{f}^q d\lambda\\
&\leq& 2^{q-1} \sum_{i \in \V} \int_{U_i} \abs{f-f_{U_i,\lambda}}^q d\lambda 
+ 2^{q-1} \sum_{i \in \V} \int_{U_i} \abs{f_{U_i,\lambda}}^q d\lambda\\
&=& 2^{q-1} \sum_{i \in \V} \int_{U_i} \abs{f-f_{U_i,\lambda}}^q d\lambda 
+ 2^{q-1} \sum_{i \in \V} \abs{f_{U_i,\lambda}}^q \lambda(U_i).
\end{eqnarray*}

The continuous Sobolev inequality gives an upper bound for the first term ; noticing that $q \geq p$ and remembering the assumptions 
on the covering, we find 
\begin{eqnarray*}
\sum_{i \in \V} \int_{U_i} \abs{f-f_{U_i,\lambda}}^q d\lambda 
&\leq&  S_c^{q/p} \sum_{i \in \V} \left( \int_{U_i^*} \abs{df}^p d\mu \right)^{q/p}\\
&\leq&  S_c^{q/p} \left( \sum_{i \in \V} \int_{U_i^*} \abs{df}^p d\mu \right)^{q/p}\\
&\leq&  S_c^{q/p} Q_1^{q/p} \left( \int_{A^\sharp} \abs{df}^p d\mu \right)^{q/p}.
\end{eqnarray*}

To estimate the second term, we use the discrete Sobolev inequality :
$$
\sum_{i \in \V} \abs{f_{U_i,\lambda}}^q \lambda(U_i)
\leq S_d \sum_{(i,j) \in \E} \abs{f_{U_i,\lambda}-f_{U_j,\lambda}}^q \max(\lambda(U_i),\lambda(U_j)).
$$
For $(i,j) \in \E$, a Hölder inequality and the fact that we have a good covering lead to :
\begin{eqnarray*}
& & \abs{f_{U_i,\lambda}-f_{U_i,\lambda}}^q \max(\lambda(U_i),\lambda(U_j)) \\
&=&  \frac{\max(\lambda(U_i),\lambda(U_j))}{\lambda (U_j)^q \lambda (U_i)^q} 
\abs{\int_{U_i} \int_{U_j} (f(x) - f(y)) d\lambda (x) d\lambda (y)}^q \\
&\leq&  \frac{\max(\lambda(U_i),\lambda(U_j))}{\lambda (U_i) \lambda (U_j)} 
\int_{U_i} \int_{U_j} \abs{f(x) - f(y)}^q d\lambda (x) d\lambda (y) \\
&\leq&  Q_2 \frac{1}{\lambda (U_{k(i,j)}^*) } 
\int_{U_{k(i,j)}^*} \int_{U_{k(i,j)}^*} \abs{f(x) - f(y)}^q d\lambda (x) d\lambda (y). \\ 
\end{eqnarray*}

Now, if $X$ is a Borel set with finite and nonzero $\lambda$-measure and if $g$ is a function in $L^q(X,\lambda)$,  
\begin{eqnarray*}
&&\frac{1}{\lambda(X)} \int_{X} \int_{X} \abs{g(x) - g(y)}^q d\lambda (x) d\lambda (y)\\ 
&\leq& \frac{1}{\lambda(X)} \int_{X} \int_{X} 2^{q-1} (\abs{g(x)}^q + \abs{g(y)}^q) d\lambda (x) d\lambda (y) \\
&\leq& 2^{q} \int_{X} \abs{g(x)}^q d\lambda (x). 
\end{eqnarray*} 

Let us apply this to $f-f_{U_{k(i,j),\lambda}^*}$, on $U_{k(i,j)}^*$ :
$$
\abs{f_{U_i,\lambda}-f_{U_j,\lambda}}^q \max(\lambda(U_i),\lambda(U_j))
\leq 
Q_2 2^{q} \int_{U_{k(i,j)}^*} \abs{f-f_{U_{k(i,j),\lambda}^*}}^q d\lambda.
$$
Now, by the continuous Sobolev inequality,
$$
\abs{f_{U_i,\lambda}-f_{U_j,\lambda}}^q \max(\lambda(U_i),\lambda(U_{j}))
\leq Q_2 2^{q} S_c^{q/p} \left( \int_{U_{k(i,j)}^\sharp} \abs{df}^p d\mu \right)^{\frac{q}{p}}. 
$$
Therefore :
$$
\sum_{i \in \V} \abs{f_{U_i,\lambda}}^{q} \lambda(U_i) 
\leq  S_d  
\sum_{(i,j) \in \E} Q_2 2^{q} S_c^{q/p} 
\left(\int_{U_{k(i,j)}^\sharp} \abs{df}^p d\mu\right)^{\frac{q}{p}}.
$$
As $q$ is greater or equal to $p$,
$$
\sum_{i \in \V} \abs{f_{U_i,\lambda}}^{q} \lambda(U_i) 
\leq  S_d Q_2 2^q S_c^{q/p} 
\left( \sum_{(i,j) \in \E} 
\int_{U_{k(i,j)}^\sharp} \abs{df}^p d\mu 
\right)^{\frac{q}{p}}.
$$
By using twice the fact that we have a good covering, we see that :
\begin{eqnarray*}
\sum_{(i,j) \in \E} 
\int_{U_{k(i,j)}^\sharp} \abs{df}^p d\mu 
&\leq& Q_1^2 \sum_{i \in \V} \int_{U_i^\sharp} \abs{df}^p d\mu  \\
&\leq& Q_1^{3} \int_{A^\sharp} \abs{df}^p d\mu.
\end{eqnarray*}
Hence :
$$
\sum_{i \in \V} \abs{f_{U_i,\lambda}}^{q} \lambda(U_i) 
\leq S_d Q_2 2^q S_c^{q/p}  Q_1^{3q/p} \left( \int_{A^\sharp} \abs{df}^p d\mu \right)^{\frac{q}{p}}.
$$

Eventually, we get :
$$
\int_A \abs{f}^q d\lambda
\leq 2^{q-1} (S_c^{q/p} Q_1^{q/p} + S_d Q_2 2^q S_c^{q/p}  Q_1^{3q/p})
\left( \int_{A^\sharp} \abs{df}^p d\mu \right)^{\frac{q}{p}}
$$
And this is what we wanted to prove.
\end{demo}

\vspace{\baselineskip}

There is also a "Neumann" version of this result.

\begin{thm}\label{recollsoboneu}
Suppose $k \in ]1,\infty]$ and $1\leq p < k$. If a finite good covering $\U$ of $A$ in $A^\sharp$ satisfies  
the continuous $L^p$ Sobolev inequality of order $k$ (\ref{sc}) and the discrete $L^p$ Sobolev-Neumann inequality of order 
$\infty$ (\ref{sdn}), the following Sobolev-Neumann inequality is true :
$$\Forall{f}{C^\infty(A)} \int_A \left( \abs{f- f_{A,\lambda}}^{\frac{p k}{k-p}} d\lambda \right)^{\frac{k-p}{k}} 
\leq S \int_{A^\sharp}  \abs{df}^p d\mu, 
$$
with 
$$
S= 2^{2p-1 + \frac{p}{k}} ((S_c Q_1)^{\frac{k}{k-p}} + S_d Q_2 (2^p S_c Q_1^3)^{\frac{k}{k-p}})^{\frac{k-p}{k}}. 
$$
\end{thm}

\begin{demo}
Again, set $q:=\frac{pk}{k-p}$ and fix $f \in C^\infty_c(A)$. 
First, note that 
$$
\norm{f - f_{A,\lambda}}_{L^q(A,\lambda)} \leq 2 \inf_{c \in \Real} \norm{f - c}_{L^q(A,\lambda)}.
$$ 
Indeed, if $c$ is a real number,
\begin{eqnarray*}
\norm{f - f_{A,\lambda}}_{L^q(A,\lambda)} &\leq& \norm{f - c}_{L^q(A,\lambda)} + \norm{c - f_{A,\lambda}}_{L^q(A,\lambda)} \\
&=& \norm{f - c}_{L^q(A,\lambda)} + \abs{f_{A,\lambda} - c} \lambda(A)^{\frac{1}{q}} \\
&=& \norm{f - c}_{L^q(A,\lambda)} + \abs{\int_A (f-c) d\lambda} \lambda(A)^{\frac{1}{q} - 1} \\
\end{eqnarray*}
By Hölder inequality,
\begin{eqnarray*}
\norm{f - f_{A,\lambda}}_{L^q(A,\lambda)} &\leq& \norm{f - c}_{L^q(A,\lambda)} + \left( \int_A \abs{f-c}^q d\lambda \right)^{\frac{1}{q}} 
\lambda(A)^{1 - \frac{1}{q}} \lambda(A)^{\frac{1}{q} - 1} \\
&=& 2 \norm{f - c}_{L^q(A,\lambda)} 
\end{eqnarray*}
As this is true for each $c \in \Real$, this proves the statement.

In particular, for 
$$
c := m_\lambda (f_{U_., \lambda}) = \frac{\sum_{i \in \V} f_{U_i,\lambda} \lambda(U_i)}{\sum_{i \in \V} \lambda(U_i)},  
$$
we can write
\begin{eqnarray*}
&\int_A& \abs{f - f_{A,\lambda}}^q d\lambda \\
&\leq& 2^q \int_A \abs{f - c}^q d\lambda \\
&\leq& \sum_{i \in \V} \int_{U_i} \abs{f-c}^q d\lambda\\
&\leq& 2^{2q-1} \sum_{i \in \V} \int_{U_i} \abs{f-f_{U_i,\lambda}}^q d\lambda 
+ 2^{2q-1} \sum_{i \in \V} \int_{U_i} \abs{f_{U_i,\lambda}-c}^q d\lambda\\
&=& 2^{2q-1} \sum_{i \in \V} \int_{U_i} \abs{f-f_{U_i,\lambda}}^q d\lambda 
+ 2^{2q-1} \sum_{i \in \V} \abs{f_{U_i,\lambda}-c}^q \lambda(U_i).
\end{eqnarray*}

We then estimate both terms as in the proof of theorem \ref{recollsobodir} : for the second, it is made possible by our choice of $c$.
\end{demo}

In fact, our argument leads to more general theorems. We will not use them but let us phrase the "Dirichlet" version. The reader 
will easily imagine the "Neumann" version. For instance, this kind of result could be used to patch local Sobolev and Poincaré 
inequalities together. 

\begin{thm}
Suppose $1\leq p \leq r \leq q \leq \infty$. Set $k=\frac{qp}{q-p}$. If a good covering $\U$ of $A$ in $A^\sharp$ satisfies 
the continuous $L^p$ Sobolev-Neumann inequality of order $k$ (with constant $S_c$),
the discrete $L^r$  Sobolev-Dirichlet inequality of order $\frac{rq}{q-r}$ (with constant $S_d$),
and the continuous $L^p$ Sobolev-Neumann inequality of order $\frac{pr}{r-p}$ (with constant $S'_c$), 
$M$ satisfies the following $L^p$ Sobolev-Dirichlet inequality of order $k$ :
$$
\Forall{f}{C^\infty_c(A)} \int_A \left( \abs{f}^{q} d\lambda \right)^{p/q} 
\leq S \int_{A^\sharp} \abs{df}^p d\mu, 
$$
with 
$$
S = 2^{p-p/q} \left( (Q_1 S_c)^{q/p} + \left( S_d Q_2 2^r (S'_c)^{r/p} \right)^{q/r} Q_1^{3q/p} \right)^{p/q} .
$$
\end{thm}


\subsection{Sobolev and isoperimetric inequalities on graphs.}

Now, we know that discrete Sobolev inequalities on appropriate graphs make it possible to patch local Sobolev inequalities together. The
problem is : how can we show such discrete inequalities ? 

Our first purpose here is to clarify the link between Sobolev inequalities of the same order on weighted graphs. We explain why, as 
in the continuous case, the $L^1$ inequality of order $k$ ($1<k\leq\infty$) imply the $L^p$ inequalities for  $1 \leq p < k$.

\begin{prop}\label{L1L2}
We consider an infinite weighted graph $(\V,\E,m)$ (see remark \ref{graphepondere}). We assume there exists  $C \geq 1$ and $d \in \Natural$ such that $$
\Forall{(i,j)}{\E} C^{-1} m(i) \leq m(j) \leq C m(i)
$$ 
and the degree of each vertex is bounded by $d$. Then the $L^1$ Sobolev inequality of order $k$, $1 < k \leq \infty$, 
\begin{equation}\label{soboL1}
\Forall{f}{L^1(\V,m)} 
\left(\sum_{i \in \V} \abs{f(i)}^{\frac{k}{k-1}} m(i)\right)^{\frac{k-1}{k}} \leq S \sum_{(i,j) \in \E} \abs{f(i)-f(j)} m(i,j).
\end{equation}
imply the $L^p$ Sobolev inequality of order $k$ for $1 \leq p < k$ :
\begin{equation}\label{soboL2}
\Forall{f}{L^p(\V,m)} \left(\sum_{i \in \V} \abs{f(i)}^{\frac{p k}{k-p}} m(i)\right)^{\frac{k-p}{p k}} 
\leq S' \left( \sum_{(i,j) \in \E} \abs{f(i)-f(j)}^p m(i,j) \right)^\frac{1}{p},
\end{equation}
where $S' = 2p \frac{k-1}{k-p} d S C^{1-\frac{1}{p}}$.
\end{prop}

\begin{preuve}
Let $f$ be an element ok $\Real^\V$ with finite support. We apply (\ref{soboL1}) to $\abs{f}^\gamma$ where $\gamma \geq 1$ is a parameter
that we will fix later :
$$
\left(\sum_{i \in \V} \abs{f(i)}^{\frac{\gamma k}{k-1}} m(i)\right)^{\frac{k-1}{k}} 
\leq S \sum_{(i,j) \in \E} \abs{\abs{f(i)}^\gamma-\abs{f(j)}^\gamma} m(i,j).
$$
If $a$, $b$ are real numbers, the following is true
$$
\abs{\abs{a}^\gamma - \abs{b}^\gamma} \leq \gamma \max(\abs{a},\abs{b})^{\gamma -1} \abs{\abs{a} - \abs{b}} 
\leq \gamma \abs{a-b} (\abs{a}^{\gamma -1} + \abs{b}^{\gamma -1}).
$$
Consequently,
$$
\left(\sum_{i \in \V} \abs{f(i)}^{\frac{\gamma k}{k-1}} m(i)\right)^{\frac{k-1}{k}} 
\leq \gamma S \sum_{(i,j) \in \E} \abs{f(i)-f(j)} (\abs{f(i)}^{\gamma-1} + \abs{f(j)}^{\gamma-1}) m(i,j).
$$
By Hölder inequality, 
\begin{eqnarray*}
& &\left(\sum_{i \in \V} \abs{f(i)}^{\frac{\gamma k}{k-1}} m(i)\right)^{\frac{k-1}{k}} \\
&\leq& \gamma S \left( \sum_{(i,j) \in \E} \abs{f(i)-f(j)}^p m(i,j) \right)^{\frac{1}{p}} 
 \Big{[} \left( \sum_{(i,j) \in \E} \abs{f(i)}^{(\gamma-1)\frac{p}{p-1}} m(i,j) \right)^{1 -\frac{1}{p}} \\
 &+& \left( \sum_{(i,j) \in \E} \abs{f(j)}^{(\gamma-1)\frac{p}{p-1}} m(i,j) \right)^{1-\frac{1}{p}} \Big{]}.
\end{eqnarray*}
And our assumptions on the graph enable us to write 
\begin{eqnarray*}
&&\left(\sum_{i \in \V} \abs{f(i)}^{\frac{\gamma k}{k-1}} m(i)\right)^{\frac{k-1}{k}} \\
&\leq& 2 \gamma d S C^{1-\frac{1}{p}} \left( \sum_{(i,j) \in \E} \abs{f(i)-f(j)}^p m(i,j) \right)^{\frac{1}{p}} 
\left( \sum_{i \in \V} \abs{f(i)}^{(\gamma-1)\frac{p}{p-1}} m(i) \right)^{1 - \frac{1}{p}}.
\end{eqnarray*}
Set $\gamma:= p \frac{k-1}{k-p} \geq 1$ to conclude the proof.
\end{preuve}

\vspace{\baselineskip}

Now, let us explain why inequalities like (\ref{soboL1}) stem from isoperimetric inequalities on the graph.

\begin{defn}
Let $(\V,\E)$ be a graph. We define the boundary $\partial \Omega$ of a subset $\Omega$ of $\V$ as  
$$
\partial \Omega := \set{(i,j) \in \E, \ \set{i,j} \cap \Omega \not= \emptyset \text{ and } 
\set{i,j} \cap (\V \backslash \Omega) \not= \emptyset}. 
$$
\end{defn}

\begin{prop}\label{isoperisobo}
Let $(\V,\E, m)$ be an infinite weighted graph and fix $k \in ]1,\infty]$. Then the isoperimetric inequality of order $k$
\begin{equation}\label{isoperi}
\forall \Omega \subset \V \text{ with } m(\Omega)<\infty, \,  m(\Omega)^{\frac{k-1}{k}} \leq I \, m(\partial \Omega).
\end{equation}
is equivalent to the $L^1$ Sobolev inequality of order $k$ 
$$
\Forall{f}{L^1(\V,m)} 
\left(\sum_{i \in \V} \abs{f(i)}^{\frac{k}{k-1}} m(i)\right)^{\frac{k-1}{k}} \leq I \sum_{(i,j) \in \E} \abs{f(i)-f(j)} m(i,j).
$$
\end{prop}

\begin{preuve}
By considering characteristic functions of subsets of $\V$, one easily sees that the Sobolev inequality implies the isoperimetric
inequality. To prove the converse, set $q=\frac{k}{k-1}$ and let $f$ be a function on $\V$, with finite support. 
For every $i \in \V$, we write
$$
f(i) = \int_0^{f(i)} dt = \int_0^\infty 1_{t < f(i)} dt.
$$
Thus,
$$
\norm{f}_{L^q(\V,m)} \leq \int_0^\infty \norm{1_{t < f(.)}}_{L^q(\V,m)} dt 
= \int_0^\infty \left( \sum_{\set{i \in \V, \ f(i) > t}} m(i) \right)^{\frac{1}{q}}  dt.
$$
If the isoperimetric inequality is true, we find
\begin{eqnarray*}
\norm{f}_{L^q(\V,m)} &\leq& I \int_0^\infty m(\partial \set{i \in \V, \ f(i) > t})  dt \\
&=& I \int_0^\infty \sum_{\set{(i,j) \in \E, \ f(j) \leq t < f(i)}} m(i,j)  dt \\
&& + I \int_0^\infty \sum_{\set{(i,j) \in \E, \ f(i) \leq t < f(j)}} m(i,j)  dt \\
&=& I \sum_{(i,j) \in \E} \abs{f(i) - f(j)} m(i,j).
\end{eqnarray*}
\end{preuve}

This paragraph shows that if the graph obtained by discretization (as explained above) satisfies an isoperimetric inequality, it will
satisfies a convenient Sobolev inequality, so that we will be able to implement our patching process.

It is time to turn to geometry so as to obtain concrete inequalities.


\section{Sobolev and Hardy inequalities on manifolds with nonnegative Ricci curvature.}

Sobolev inequalities are a major tool of global analysis. Unfortunately, they are not always available. It is known that on manifolds with
nonnegative Ricci curvature and maximal volume growth, they actually occur (\cite{Cro}), providing a lot of analytical, 
geometrical and topological information : see \cite{BKN}, for instance. As soon as the volume growth is not maximal, the Sobolev 
inequality cannot be true. Our aim here is to show that, even if the volume growth is not maximal, a weighted Sobolev inequality still occurs. 

\subsection{Geometric preliminaries.}
 
We would like to emphasize here some features of complete manifolds with nonnegative Ricci curvature. These are the typical manifolds where
our discretization scheme applies. 

Recall that if $x$ is a point in $M$, we denote by $V(x,R)$ the volume of the ball $B(x,R)$ centered in $x$ and with radius $R$. 
We will sometimes omit the center when it is some distinguished point $o$. We also introduce $A(R_1,R_2) := B(R_2) \backslash B(R_1)$ and 
$V(R_1,R_2) := \vol A(R_1,R_2)$.

First, the Bishop-Gromov comparison theorem says that, in manifolds with nonegative Ricci curvature, the volume growth of balls 
is "subeuclidian" in a very strong way.

\begin{thm}[Bishop-Gromov]
Let $M$ be a complete riemannian manifold with nonnegative Ricci curvature. Then for every $x$ in $M$, the function $\rho_x$ defined for $t
\geq 0$ by
$$
\rho_x(t) = \frac{t^n}{V(x,t)}
$$ 
is a nondecreasing function. It implies that for $0<r<t$, 
\begin{equation}\label{bishopgromov}
\Forall{x}{M} \frac{\vol B(x,t)}{\vol B(x,r)} \leq \left( \frac{t}{r} \right)^n.
\end{equation}
And a useful corollary is that for $x,y \in M$ and $0<r<t+d(x,y)$ :
\begin{equation}\label{bishopgromovdecale}
\frac{\vol B(y,t)}{\vol B(x,r)} \leq \frac{\vol B(x,t+d(x,y))}{\vol B(x,r)} \leq \left( \frac{t+d(x,y)}{r} \right)^n.
\end{equation}
\end{thm}

For a proof, see \cite{Cha}.

Note the following simple consequence. The argument of the proof will constantly be used in the sequel.
\begin{cor}\label{estimanneaux1}
Let $M^n$ be a connected complete noncompact riemannian manifold with nonnegative Ricci curvature. 
Then for every $\kappa>1$, there exists a positive constant $C(n,\kappa)$ such that for any $x \in M$ and $R> 0$,
$$
C(n,\kappa)^{-1} 
\leq \frac{\vol B(x,\kappa R) \backslash B(x,R) }{\vol B(x,R) \backslash B(x,\kappa^{-1}R) } 
\leq C(n,\kappa).
$$
\end{cor}

\begin{preuve}
To prove the lower bound, choose a point $y$ on the sphere $S(x,(\kappa +1)R/2)$ centered in $x$ and of radius 
$(\kappa+1)R/2$ (such a point exists since $M$ is assumed to be complete, noncompact and connected). Then the 
ball $B:=B(y,(\kappa - 1)R/2)$ is contained in $B(x,\kappa R) \backslash B(x,R)$. Therefore
$$
\frac{\vol( B(x,R) \backslash B(x,\kappa^{-1}R) )}{\vol( B(x,\kappa R) \backslash B(x,R) )} 
\leq \frac{\vol B(x,R)}{\vol B(y,(\kappa - 1)R/2)},
$$
and (\ref{bishopgromovdecale}) yields
$$
\frac{\vol( B(x,R) \backslash B(x,\kappa^{-1}R) )}{\vol( B(x,\kappa R) \backslash B(x,R) )} 
\leq \left( \frac{R + \frac{(\kappa +1)R}{2}}{\frac{(\kappa - 1)R}{2}}\right)^n 
=\left( \frac{\kappa + 3}{\kappa- 1} \right)^n.
$$
The upper bound is proved likewise.
\end{preuve}

\vspace{\baselineskip}

Moreover, starting from the comparison theorem, P. Buser \cite{Bus} showed the following  
\begin{thm}[Buser]
In a complete noncompact riemannian manifold with nonnegative Ricci curvature, for any $p$ in $[1,\infty[$, every ball $B(x,R)$ satisfies 
the $L^p$ Poincaré inequality
\begin{equation}\label{buser}
\Forall{f}{C^\infty(B(x,R))} \int_{B(x,R)} \abs{f - f_{B(x,R)}}^p dvol  
\leq C(n,p) R^p \int_{B(x,R)} \abs{df}^p dvol,
\end{equation}
where $f_{B(x,R)}$ denotes the mean value of $f$ on the ball $B(x,R)$, with respect to the riemannian measure $vol$.
\end{thm}

This result yields the fundamental inequalities we need. Besides, it will prove useful in the study of the geometry at infinity 
of manifolds with nonnegative Ricci curvature.

Let us mention the Cheeger-Gromoll theorem (\cite{CG},\cite{Bes}), which enlightens the structure of manifolds with 
nonnegative Ricci curvature :
\begin{thm}[Cheeger-Gromoll]
A connected complete riemannian manifold with nonnegative Ricci curvature is always the riemannian product of the euclidian space $\Real^d$  
and a connected complete riemannian manifold with nonnegative Ricci curvature which possesses no line.
\end{thm}

\begin{cor}
A connected complete noncompact riemannian manifold with nonnegative Ricci curvature possesses exactly one end, 
unless it is a riemannian product of $\Real$ and a compact manifold.
\end{cor}

\begin{rem}
In our setting, the volume growth of balls will forbid the particular case, which is therefore irrelevant here.
\end{rem}

In what follows, we will be working on annuli so that we are interested in their
topology/geometry, and especially in their connectedness : it is an obvious necessary condition if we hope to show a Sobolev or Poincaré inequality on them. In \cite{And}, M. Anderson proved that the first Betti number of a connected complete riemannian manifold with nonnegative Ricci curvature is bounded by its dimension. Now, \cite{LT} points out a consequence of the finiteness of the first Betti number :
\begin{prop}\label{litam}
Let $M$ be a connected complete riemannian manifold with nonnegative Ricci curvature, finite first Betti number
and exactly $k$ ends. Let us fix a point $o \in M$ and consider balls and annuli centered in $o$. Then for large $R$ and any $r>0$, 
denoting by $M_R$ the union of all unbounded connected components $M \backslash B(R)$, it is true that 
$A(R,R+r) \cap M_R$ has exactly $k$ connected components.
In particular, if $M$ has exactly one end, for large $R$ and any $r>0$, the annulus $A(R,R+r)$ possesses one and only one
component that can be connected to infinity inside $M \backslash B(R)$.
\end{prop}

Let us give an interpretation in terms of discretization. Consider a manifold $M$ with nonnegative Ricci curvature, possessing one end, and fix a point in $M$. Let us choose $R>0$ and $\kappa>0$. We discretize $M$ in the following manner. We associate a vertex to $B(R)$ and to every connected component of the annuli $A(\kappa^i R,\kappa^{i+1} R)$, $i \in \Natural$. Let us decide that there is an edge between two given vertices if and only if the closures of the corresponding subsets of $M$ intersect. Then the proposition above says that for large $R$ this graph is a tree and its root is the vertex corresponding to $B(R)$. From another point of view, it says, that even if $R$ is small, outside a finite subset, the graph is a tree.

Now, there is no reason why this tree should not have branches, and for technical reasons (see the proof of lemma \ref{sobodisc} below), we would like to make them as small as possible. What we need is some kind of control on the size of bounded connected components of the complements of balls in the manifold. This is given by the following proposition, which we state with rather general assumptions. 

\begin{prop}[RCA]
Let $M$ be a connected complete riemannian manifold, satisfying the doubling volume property
$$
\Forall{x}{M} \forall \, R>0, \,  V(x,2R) \leq C_D V(x,R),
$$
the scaled $L^p$ Poincaré inequality centered in some point $o$ in $M$ 
$$
\Forall{f}{C^\infty_c(M)} \forall \, R>0, \,  \int_{B(o,R)} \abs{f - f_{B(o,R)}}^p dvol \leq C_P R^p \int_{B(o,R)} \abs{df}^p dvol
$$
and the inverse doubling volume condition centered in $o$
$$
\forall R_2 \geq R_1 > 0, \  \frac{V(o,R_2)}{V(o,R_1)} \geq C_o \left( \frac{R_2}{R_1} \right)^\nu
$$
with $\nu > p$. Here, $C_D \geq 1$, $p \geq 1$, $C_P >0$, $C_o >0$. Then there exists $\kappa_0 >0$ such that for $R > 0$, 
if $x,y$ are two points in $S(o,R)$, there is a path from $x$ to $y$ which remains inside $B(o,R)\backslash B(o,\kappa_0^{-1} R)$. Moreover,
it is possible to find an explicit constant, in terms of $p,C_D,C_P,C_o,\nu$.  
\end{prop}

In terms of the discretization we have introduced, this means that for large $\kappa$, for every two vertices on the same level of the tree 
(i.e. corresponding to the same annulus), there exists a vertex of the previous level which is connected to both of them. 

\noindent
\begin{figure}
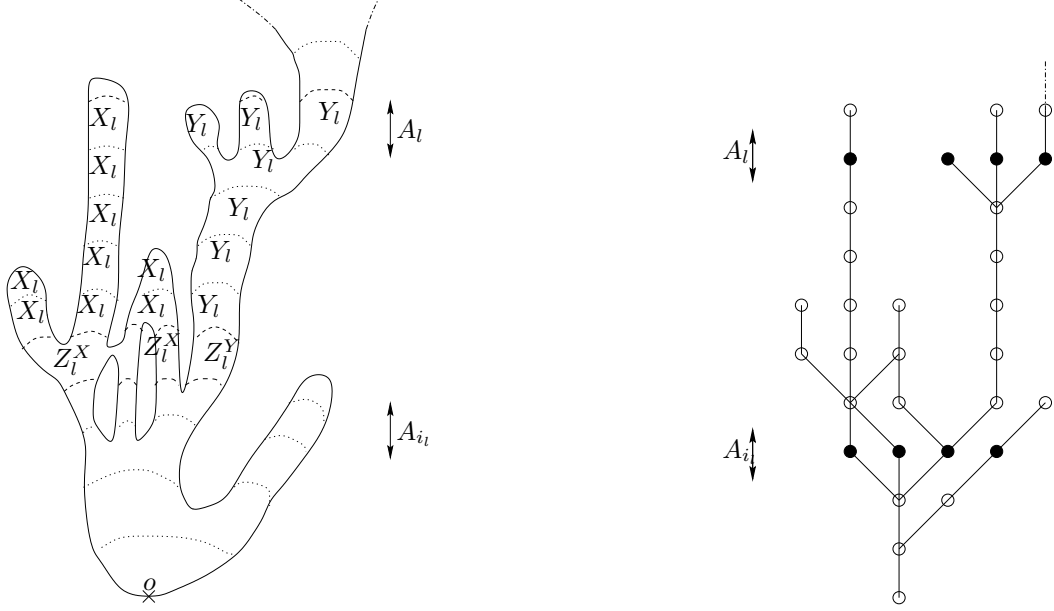

\psfrag{X}{$X_l$}  \psfrag{Y}{$Y_l$}  \psfrag{ZX}{$Z_l^X$}  \psfrag{ZY}{$Z_l^Y$}  \psfrag{T}{ }  \psfrag{Al}{$A_l$}  \psfrag{Ai}{$A_{i_l}$}  \psfrag{o}{$o$}
\includegraphics[width=5.4cm]{cactus2.pstex}
\hspace{4cm}
\psfrag{S1}{ }  \psfrag{S2}{ }  \psfrag{S3}{ }  \psfrag{S4}{ }  \psfrag{T}{ }
\includegraphics[width=4.4cm]{arbre2.pstex}
\label{fig1}
\caption{A manifold and its discretization.} 
\end{figure}

\begin{preuve}
Consider the graph obtained as above by working with $A_i := A(2^{i-1} R,2^{i} R)$, $i \in \Natural^*$, $R>0$, plus $B(R)=:A_{0}$. Set 
$B_i = B(2^{i} R)$. We define $\C$ as the bijective map which associates to every vertex of the graph the corresponding component of annulus. 
Let us write $\A_i$ for $\C^{-1}(A_i)$ and $\B_i$ for $\C^{-1}(B_i)$.

Now, fix $l \in \Natural^*$. We consider the nonempty set
$$
I_l = \set{i \in [0,l], \, \A_l \text{ is contained in a connected component of } \B_l \backslash \B_{i-1}}
$$
and set $i_l = \max I_l$. Call $\M_l$ the connected component of $\B_l \backslash \B_{i_l-1}$ which contains $\A_l$. We assume $l-i_l$ is
greater than $3$ and think of it as a large number.

By definition, $\M_l \backslash \A_{i_l}$ is not connected : we choose one of its connected component $\X'_l$ and name $\Y'_l$ the union of the
other connected components. We finally define $X'_l :=\C^{-1} (\X'_l)$, $Y'_l :=\C^{-1} (\Y'_l)$, $X_l := X'_l \backslash A_{i_l +1}$, 
$Y_l := Y'_l \backslash A_{i_l +1}$, $Z^X_l := X'_l \cap A_{i_l +1}$,
$Z^Y_l := Y'_l \cap A_{i_l +1}$ and $Z_l := Z^X_l \cup Z^Y_l$ (see figure \ref{fig1}). 

Given real numbers $a$ and $b$, we can define a Lipschitz function $f_l$ on $B_l$ in the following way : 
$$
f_l = \left\{ \begin{array}{lllll}
a  & \text{on $X_l$,} \\
b  & \text{on $Y_l$,} \\
a \frac{r_o - 2^{i_l} R}{2^{i_l} R} & \text{on $Z^X_l$,} \\
b \frac{r_o - 2^{i_l} R}{2^{i_l} R} & \text{on $Z^Y_l$,} \\
0 & \text{everywhere else.}
\end{array}\right.
$$

The Poincaré inequality says
\begin{equation}\label{bus}
\int_{B_l} \abs{f_l - (f_l)_{B_l}}^p dvol  
\leq C_P 2^{lp} R^p \int_{B_l} \abs{df_l}^p dvol.
\end{equation}
We choose $a$ and $b$ so that the mean value of $f_l$ on $X_l \cup Y_l$ is $0$ : 
$$
a \vol X_l + b \vol Y_l = 0.
$$
With $a:=1$, this means $b = -\frac{\vol X_l}{\vol Y_l}$.

On the one hand,
\begin{eqnarray*}
\int_{B_l} \abs{f_l - (f_l)_{B_l}}^p dvol 
&\geq& 2^{-p} \frac{\int_{B_l} \int_{B_l} \abs{f_l(x) - f_l(y)}^p dx dy}{\vol B_l}  \\
&\geq& 2^{-p} \frac{\vol X_l \vol Y_l \abs{b-a}^p}{\vol B_l}  \\
&=& 2^{-p} \frac{\vol X_l \vol Y_l \left( 1 + \frac{\vol X_l}{\vol Y_l} \right)^p}{\vol B_l}.
\end{eqnarray*}

On the other hand,
\begin{eqnarray*}
\int_{B_l} \abs{df_l}^p dvol 
&\leq& \left( \vol Z^X_l \left( \frac{a}{2^{i_l} R} \right)^p + \vol Z^Y_l \left( \frac{b}{2^{i_l} R} \right)^p \right)\\
&=&  \frac{\vol Z^X_l + \vol Z^Y_l \left( \frac{\vol X_l}{\vol Y_l} \right)^p }{2^{i_lp} R^p}.
\end{eqnarray*}

So 
\begin{eqnarray*}
2^{-p} \frac{\vol X_l \vol Y_l \left( 1 + \frac{\vol X_l}{\vol Y_l} \right)^p}{\vol B_l} 
&\leq& C_P 2^{p(l-i_l)} \left( \vol Z^X_l + \vol Z^Y_l \left( \frac{\vol X_l}{\vol Y_l} \right)^p \right)\\
&\leq& C_P 2^{p(l-i_l)} \vol Z_l \left( 1 + \left( \frac{\vol X_l}{\vol Y_l} \right)^p \right),
\end{eqnarray*}
hence 
\begin{equation}\label{ineg}
1 \leq 2^p C_P 2^{p(l-i_l)} \frac{\vol Z_l \vol B_l}{\vol X_l \vol Y_l}. 
\end{equation}

Now
$$
\vol Z_l \leq V(o,2^{i_l+1} R).
$$

A lower bound on $\vol X_l$ can be obtained as in the proof of (\ref{estimanneaux1}). Choose a point $x_l$ in
$S(o,(2^{l-2} + 2^{l-1})R/2) \cap X_l$ et note that $B(x_l,2^{l-3}R)$ is contained in $X_l$ : it lies in $A(2^{l-2}R, 2^{l-1}R)$ 
and it is connected, so it lies in the connected component of its center $x_l$ in $A(2^{l-2}R, 2^{l-1}R)$, hence in $X_l$.
The doubling volume property implies
$$
\Forall{x}{M} \forall \, R_2 \geq R_1>0, \,  V(x,R_2) \leq C_D (R_2/R_1)^{\log_2 C_D} V(x,R_1),
$$
so that
$$
\frac{V(o,2^{l}R)}{V(x_l,2^{l-3}R)} 
\leq C_D \left( \frac{ 2^{l} + (2^{l-2} + 2^{l-1})/2}{2^{l-3} } \right)^{\log_2 C_D}
= C_D 11^{\log_2 C_D}
$$
and 
$$
\vol X_l \geq V(x_l,2^{l-3}R)
\geq C_D^{-1} 11^{-\log_2 C_D} V(o,2^{l}R).
$$ 

As we have the same lower bound for $\vol Y_l$, (\ref{ineg}) gives :
$$
1 \leq 2^p C_P C_D^2 121^{\log_2 C_D} 2^{p(l-i_l)} \frac{V(o,2^{i_l+1}R)}{V(o,2^l R)}.
$$
(\ref{volboules1}) enables us to write :
$$
1 \leq 2^p C_P C_D^2 121^{\log_2 C_D} 2^\nu C_o 2^{(l-i_l)(p-\nu)}.
$$
Since $\nu > p$, this inequality says that $l-i_l$ is bounded by some constant independent of $l$ : 
the branches of the tree have a bounded length. (\ref{RCA}) stems from it easily.
\end{preuve}

\vspace{\baselineskip}

\begin{cor}\label{RCA}
Let $M$ be a connected complete riemannian manifold with nonnegative Ricci curvature and assume there are $o$ in $M$, $C_o>0$ and $\nu > 1$
such that
$$
\forall R_2 > R_1 \geq 1, \  \frac{V(o,R_2)}{V(o,R_1)} \geq C_o \left( \frac{R_2}{R_1} \right)^\nu
$$
Then there exists $\kappa_0 = \kappa_0(n,\nu, C_o) >0$ such that for $R > 0$, 
if $x,y$ are two points in $S(o,R)$, there is a path from $x$ to $y$ which remains inside $B(o,R)\backslash B(o,\kappa_0^{-1} R)$. 
\end{cor}


\subsection{Poincaré and Sobolev inequalities on connected components of annuli.}

We show here that Poincaré or Sobolev inequalities on balls imply analogous inequalities on connected subsets of annuli.

We first use Buser theorem (\ref{buser}) :

\begin{lem}\label{poincann}
Let $M^n$ be a noncompact connected complete riemannian manifold with nonnegative Ricci curvature.
Fix $p\geq 1$, $R >0$, $\kappa>1$ and consider a connected Borel subset $A$ of the annulus $B(o,\kappa R) \backslash B(o,R)$, $o \in M$. Then if 
we let $A_\delta$ be the $\delta R$-neighbourhood of $A$, with $0<\delta<1$, the following Poincaré inequality is true :
$$
\Forall{f}{C^\infty (A_\delta)} \int_{A} \abs{f - f_{A}}^p dvol  
\leq C(n,\kappa,\delta,p) R^p \int_{A_\delta} \abs{df}^p dvol.
$$
\end{lem}

\begin{preuve}
Set $s = \delta R$ and consider a $s$-lattice $(x_i)_{i \in I}$ of $A$, i.e. a maximal subset of 
$A$ such that the distance between any two of its elements is at least $s$. We then set $V_i = B(x_i,s)$, 
$V_i^* = V_i^\sharp= B(x_i,3s)$.

It is easy to see that $(V_i,V_i^*,V_i^\sharp)_{i \in I}$ is a good covering of $A$ in $A_\delta$ (cf. (\ref{bonrec})), 
with respect to the riemannian measure. Indeed, for (iii), we can note that the $V_i^*$ under consideration are contained in $B(x_{i_0},9s)$
and use (\ref{bishopgromovdecale}) to get $\vol(B(x_{i_0},9s)) \leq 30^n \vol(B(x_{i},\frac{s}{2}))$ ; since the balls  $B(x_{i},\frac{s}{2})$ 
do not intersect,  we see that $Q_1 = 30^n$ is convenient. In (iv), we can choose $k(i,j) = i$. As to (v), (\ref{bishopgromov}) yields 
$\vol(V_{i}^*) \leq 3^n \vol(V_{i})$ and (\ref{bishopgromovdecale}) gives $\vol(V_{i}^*) \leq 5^n \vol(V_{j})$, so that we can set $Q_2 =
5^n$.

We intend to apply the theorem \ref{recollsoboneu} with $k=\infty$. First, (\ref{buser}) yields the continuous   
inequality, with constant $C(n,p) s^2$. What about the discrete inequality ? 

Noticing the balls  $B(x_i,\frac{s}{2})$ do not intersect and are contained in the ball $B(o,\kappa R + \frac{s}{2} )$, we find that
$$
\Card(I) \min_{i \in I} \vol(B(x_i,s/2)) \leq \vol(B(o,\kappa R + s/2)),
$$
and with (\ref{bishopgromovdecale}), this implies an upper bound on the number of balls in the covering
$$
\Card(I) \leq \left( \frac{\kappa R + s/2 + \kappa R}{s/2} \right)^n = 
(1+4\kappa/\delta)^n =: N = N(n,\kappa,\delta).
$$
The point is it is independent of $R$.

Now, every finite connected graph endowed with the counting measure satisfies a Poincaré inequality : this stems from the fact that any two
norms on a vector space of finite dimension are equivalent (the connectivity is necessary here to ensure that we indeed compare two norms).
As there is only a finite number of such graphs which have at most $N$ vertices, we conclude that every such graph satisfies a Poincaré
inequality for some constant $P=P(N,p)$ (see below for an explicit constant). 
Since (\ref{bishopgromovdecale}) implies 
$$
\Forall{i,j}{\V} \frac{\vol(V_i)}{\vol(V_j)} \leq (1+2\kappa/\delta)^n,
$$
there is a number $K=K(n,\kappa,\delta)\geq 1$ such that
$$
K^{-1} m_0 \leq m(i) \leq K m_0, 
$$
where $m_0$ is proportionnal to the counting measure on our graph $\G =(\V,\E)$.

Then for every $f \in \Real^{\V}$ :
\begin{eqnarray*}
\left(\sum_{i \in \V} \abs{f(i) - m(f)}^p m(i) \right)^{1/p}
 & \leq & 2 \inf_{c \in \Real} \left(\sum_{i \in \V} \abs{f(i) - c}^p m(i) \right)^{1/p} \\
 & \leq & 2 \left(\sum_{i \in \V} \abs{f(i) - m_0(f)}^p m(i) \right)^{1/p} \\
 & \leq & 2 K \left(\sum_{i \in \V} \abs{f(i) - m_0(f)}^p m_0(i) \right)^{1/p} \\
 & \leq & 2 P K^{1/p} \left(\sum_{(i,j) \in \E} \abs{f(i)-f(j)}^p m_0(i,j) \right)^{1/p} \\
 & \leq & 2 P K^{2/p} \left(\sum_{(i,j) \in \E} \abs{f(i)-f(j)}^p m(i,j) \right)^{1/p}
\end{eqnarray*}

This yields a discrete Poincaré inequality with a constant depending only on $n,\kappa,\delta,p$ and finishes the proof, 
thanks to (\ref{recollsoboneu}).
\end{preuve}

\vspace{\baselineskip}

The same pattern gives an analogous Sobolev inequality. We first recall a theorem of L. Saloff-Coste (\cite{SC1}, \cite{SC2}) : in a
complete riemannian manifold with nonnegative Ricci curvature, every smooth function $f$ on a ball $B(x,R)$ satisfies the Sobolev inequality
\begin{equation}\label{saloff}
\left( \int_{B(x,R)} \abs{f - f_{B(x,R)}}^{\frac{2n}{n-2}} dvol \right)^{\frac{n-2}{n}}  
\leq C(n) \frac{R^2}{V(x,R)^{\frac{2}{n}}} \int_{B(x,R)} \abs{df}^2 dvol,
\end{equation}

Note that this result follows in fact from (\ref{buser}) and (\ref{bishopgromov}). We deduce the 

\begin{lem}\label{soboann}
Let $M^n$ be a noncompact connected complete riemannian manifold with nonnegative Ricci curvature and $n \geq 3$. Fix $R >0$, 
$\kappa>1$ and consider a connected Borel subset $A$ of the annulus $B(o,\kappa R) \backslash B(o,R)$, $o \in M$. Then if we let $A_\delta$ be the $\delta R$-neighbourhood of $A$, with $0<\delta<1$, the following Sobolev inequality is true.
$$
\Forall{f}{C^\infty (A_\delta)}
\left(\int_{A} \abs{f - f_{A}}^{\frac{2n}{n-2}} dvol \right)^{\frac{n-2}{n}} 
\leq C(n,\kappa,\delta) \frac{R^2}{V(o,R)^{2/n}} \int_{A_\delta} \abs{df}^2 dvol.
$$
\end{lem}

\begin{preuve}
We just explain how to adapt the previous argument, using the same notation. We set $q=\frac{2n}{n-2}$.

We want to apply (\ref{recollsoboneu}) for $p=2$ and $k=n$, with the same good covering. The discrete $L^q$ Poincaré inequality we need 
is given by the previous proof. 

(\ref{bishopgromovdecale}) gives for every $i$ in $I$
$$
\frac{V(o,R)}{V(x_i,\delta R)} 
\leq \left( \frac{1 + \kappa}{\delta} \right)^n,
$$
hence $V(x_i, 3s) \geq V(x_i,s) \geq C(n,\kappa,\delta) V(o,R)$, so that the Saloff-Coste theorem (\ref{saloff}) 
yields a continuous Sobolev-Neumann inequality for the pair of measures $(\vol,\vol)$ :
\begin{equation}
\Forall{f}{C^\infty (V_i^*)}  \left(\int_{V_i} \abs{f - f_{V_i}}^{q} dvol \right)^{\frac{2}{q}} 
\leq C(n,\kappa,\delta) R^2 V(o,R)^{-2/n} \int_{V_i^*} \abs{df}^2 dvol,
\end{equation}  
and 
\begin{equation}
\Forall{f}{C^\infty(V_i^\sharp)}  \left(\int_{V_i^*} \abs{f - f_{V_i^*}}^{q} dvol \right)^{\frac{2}{q}} 
\leq C(n,\kappa,\delta) R^2 V(o,R)^{-2/n} \int_{V_i^\sharp} \abs{df}^2 dvol.
\end{equation}  
(\ref{recollsoboneu}) ends the proof.
\end{preuve}

\vspace{\baselineskip}

Let us make a little remark. In the arguments above, we claimed the existence of the constants $P$ and $S$. Indeed, we can make them explicit,
using the following proposition.

\begin{prop}
Consider a finite connected graph $\G = (\V,\E)$ with $N_v$ vertices, endowed with the counting measure. Fix $p \geq 1$. Then 
for every real function $f$ on $\V$,
$$
\sup_{i \in \V} \abs{f(i) - m(f)} \leq N_e^{1-1/p}  \left( \sum_{(i,j) \in \E} \abs{f(i) - f(j)}^p \right)^{1/p} 
$$
and in particular,
$$
\sum_{i \in \V} \abs{f(i) - m(f)}^p \leq N_v (N_v-1)^{p-1} \sum_{(i,j) \in \E} \abs{f(i) - f(j)}^p. 
$$ 
\end{prop}

\begin{preuve}
First, we can assume the graph is a tree : cutting off edges does not change the left-hand sides and makes the right-hand sides of the
inequalities grow. Then, the graph has exactly $N_v-1$ edges. Now, to each edge $e$ we associate a copy $I_e$ of the segment $[0,1]$ ; the 
ends of $I_e$ (corresponding to $0$ and $1$) 
can be viewed as two vertices in the graph $\G$. We then build a space $X$ by gluing all $I_e$, in the natural
way, that is, we decide that the ends of segments corresponding to the same vertex in $\G$ give rise to one point in $X$. $X$ is endowed with
a natural topology and a natural Borel measure, steming from those of $[0,1]$. Note that the complement $\tilde{X}$ of the
points where two segments are glued together even possesses a natural differential structure, and a riemannian metric. 
Given $f \in \Real^\V$, we can define a continuous function $g$ on $X$ in the following manner : $g$ is linear on each segment $I_e$ and its
values at the ends of segments are simply those of $f$. Let $e$ be an edge of the graph, between the vertices $i$ and $j$, that we identify
(respectively) with $0$ and $1$ in $[0,1]$. The restriction of $g$ on $I_e$ can be identified with a function $g_e$ defined on $[0,1]$ 
by the formula 
$$
g_e(t) = f(i) + t (f(j) - f(i)).
$$
Such a function $g$ has a derivative $g'$ which is defined outside the vertices and constant on the (image in $X$ of the) interior of each 
$I_e$ : $g'_e$ = f(j) - f(i).
We claim the following inequality is true
\begin{equation}\label{poincgraph}
\norm{g}_{L^\infty(X)} \leq (N_v-1)^{1-1/p}  \norm{g'}_{L^p(X)}
\end{equation}
as soon as $g$ is continuous on $X$, $C^1$ on $\tilde{X}$ and vanishes somewhere. Let us prove it. We choose $x_0$ such that $g(x_0)=0$. Then,
given a point $x$ in the arcwise connected space $X$, we can find a unit speed path $\gamma$ from $x_0$ to $x$ which runs along each segment 
at most once. We can write
$$
g(x) = \int_\gamma g'
$$   
and the Hölder inequality implies
$$
\abs{g(x)} \leq \length(\gamma)^{1-1/p} \left( \int_\gamma \abs{g'}^p \right)^{1/p} 
\leq (N_v-1)^{1-1/p}  \norm{g'}_{L^p(X)}.
$$
Now, we want to apply this to the function $g \in C^0(X)$ which is obtained from a function $f \in \Real^\V$ with zero mean
value. As $g$ takes every value in the convex hull of the values of $f$, such a $g$ vanishes at some point, so that $g$ satisfies 
(\ref{poincgraph}). Eventually, we observe
$$
\norm{g}_{L^\infty(X)} = \norm{f}_{L^\infty(\V)}
$$ 
and
$$
\norm{g'}_{L^p(X)} = \left( \sum_{(i,j) \in \E} \abs{f(i) - f(j)}^p \right)^{1/p},
$$
and we are done.
\end{preuve}

\begin{rem}
It is possible to give a discrete proof of this result. For instance, observing that for any real number $c$
$$
\left( \sum_{i \in \V} \abs{f(i) - m(f)}^p \right)^{1/p} \leq 2 \left( \sum_{i \in \V} \abs{f(i) - c}^p \right)^{1/p}
$$
we can choose $c$ so that $f - c$ vanishes at some vertex. It is then easy to 
adapt the argument above, keeping it completely discrete. But the constant we find that way is twice the one in the proposition.
\end{rem}


\subsection{The weighted Sobolev inequality.}

In this paragraph, $M$ is a connected complete riemannian manifold, with dimension $n\geq 3$, nonnegative Ricci curvature
and satisfying (\ref{volboules1}) for some point $o$. We want to prove a weighted Sobolev inequality on $M$, by applying the theorem (\ref{recollsobodir}) for $p=2$ and $k=n$ with a good covering that we design now. 

\subsubsection{A good covering}\label{goodcov}

We fix some large $\kappa$, so as to be
sure that, for any $R > 0$, any two connected components of $A(R,\kappa R)$ are contained in one connected component of 
$A(\kappa^{-1}R,\kappa R)$ : this is made possible by (\ref{RCA}). Recall $\kappa$ can be chosen so that it depends only on 
$n$, $C_o$ and $\nu$. We also choose a ray starting from $o$ and call it $\gamma$. We will sometimes use the notation 
$R_i := \kappa^{i}$, $i \in \Integer$.

For every integer $i$, we denote by $U'_{i,a}$, $0 \leq a \leq h'_i$ the connected components of $A(R_{i-1}, R_i)$, $U'_{i,0}$ being 
the one which intersects $\gamma$. As in the proofs of \ref{poincann} and \ref{soboann}, (\ref{bishopgromovdecale}) provides a bound 
$h=h(n,\kappa)<\infty$ on the various $h'_i$, $i \in \Integer$.

A priori, this will not yield a good covering because some of the $U'_{i,a}$ may be small compared to their neighbours, contradicting (v) 
in \ref{bonrec}. This is the reason why we need to modify the covering slightly : we will glue every small component on the
level $i$ to a large one on the level $i-1$. Let us explain what we mean precisely.

We proceed in two steps. 
\begin{itemize}
\item First, we set $U_{i,a} = U'_{i,a}$ for every $i \in \Integer$ and  $1 \leq a \leq h'_i$ such that 
$\overline{U'_{i,a}}$ intersects $A(R_i,R_{i+1})$ ; every such $U_{i,a}$ contains a point $x$ on the 
sphere $S((R_{i-1}+R_i)/2)$ and thus a ball centered in $x$ and with radius $R_{i-2}$, whose volume is comparable to $V(R_i)$ (with 
(\ref{bishopgromovdecale})). 
\item Then we consider every $(i,a)$ such that $\overline{U'_{i,a}} \cap A(R_i,R_{i+1})$ is empty. 
There is $b$ in $[0,h'_{i-1}]$ such that $U'_{i,a} \cup U'_{i-1,b}$ is connected : we enlarge $U_{i-1,b}$ by adding $U'_{i,a}$ to it. 
\end{itemize}

After deleting the indices which are not used any more, this yields a covering 
$(U_{i,a})$ of $M\backslash \set{o}$, indexed by $i \in \Integer$ and $a \in [0,h_i]$, $h_i \leq h'_i$, with 
$U_{i,a} \subset A(R_{i-1},R_{i+1})$ and $\vol U_{i,a} \approx V(R_i)$. 

The following figure gives an example : on the left, different connected components of annuli $A(R_{i-1}, R_i)$ ; in the 
center, the modified covering ; on the right, the associated graph.  
\noindent\begin{center}
\centering  \psfrag{U10}{$U'_{1,0}$}\psfrag{U20}{$U'_{2,0}$}\psfrag{U30}{$U'_{3,0}$}\psfrag{U40}{$U'_{4,0}$}\psfrag{U31}{$U'_{3,1}$}\psfrag{U32}{$U'_{3,2}$}\psfrag{U41}{$U'_{4,1}$}
\psfrag{V10}{$U_{1,0}$}\psfrag{V20}{$U_{2,0}$}\psfrag{V30}{$U_{3,0}$}\psfrag{V40}{$U_{4,0}$}\psfrag{V31}{$U_{3,1}$}
\psfrag{10}{$(1,0)$}\psfrag{20}{$(2,0)$}\psfrag{30}{$(3,0)$}\psfrag{40}{$(4,0)$}\psfrag{31}{$(3,1)$}
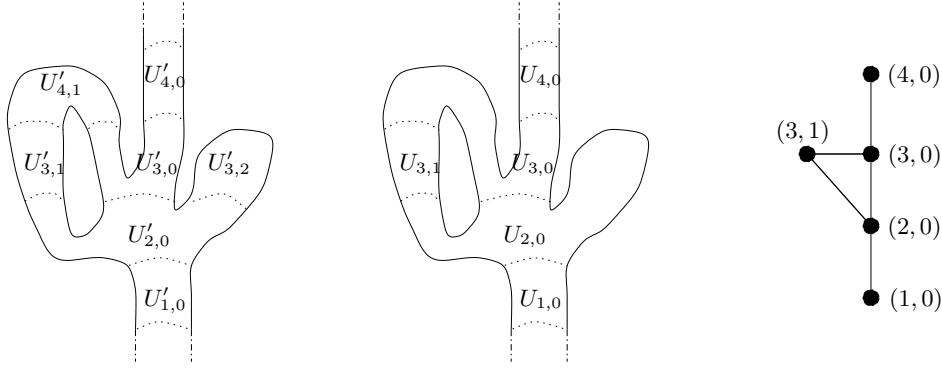
\includegraphics[width=12cm]{cactuszoom.pstex}
\label{fig2}
\label{example}
\end{center}

For $i$ in $\Integer$ and $0 \leq a \leq h_i$, we furthermore let $U_{i,a}^*$ be the union of all the $U_{j,b}$, $j \in \Integer$, $0 \leq b \leq h_j$, whose closure intersects $\overline{U_{i,a}}$. And likewise, let $U_{i,a}^\sharp$ be the union of all the $U_{j,b}^*$, $j \in \Integer$, $0 \leq b \leq h_j$, whose closure intersects $\overline{U_{i,a}^*}$. 

Now, we introduce the measure
$$
d\mu_\rho = \rho(r)^{-\frac{2}{n-2}} dvol,
$$
where $r=r_o=d(o,.)$ and $\rho(t)=\rho_o(t)$ is defined for $t \geq 0$ by  
$$
\rho(t) = \frac{t^n}{V(t)}.
$$ 
Bishop-Gromov theorem says it is a nondecreasing function and indeed, for $0<R_1 \leq R_2$,
\begin{equation}\label{controle}
1 \leq \frac{\rho(R_2)}{\rho(R_1)} \leq \left( \frac{R_2}{R_1} \right)^n ;
\end{equation}
besides, $\rho(0) = \omega_n^{-1}$, where $\omega_n$ denotes the volume of the unit
sphere in $\Real^n$. 

It is easy to see that $\U=(U_{i},U_{i}^*,U_{i}^\sharp)$ is a good covering of $M$ in $M$ with respect to $(\mu_\rho,\vol)$ : (v) is again a consequence of (\ref{bishopgromovdecale}). 

Let us prove the continuous and discrete Sobolev inequalities we need.


\subsubsection{The continuous Sobolev inequality.}

\begin{lem}\label{sobocont}
For every $i \in \Integer$ and $0 \leq a \leq h_i$, each smooth function $f$ on 
$U_{i,a}^\sharp$ satisfies
$$
\left( \int_{U_{i,a}} \abs{f - f_{U_{i,a}}}^{\frac{2n}{n-2}} d\mu_\rho \right)^{\frac{n-2}{n}} 
\leq S_c \int_{U_{i,a}^*} \abs{df}^2 dvol
$$
and
$$
\left( \int_{U_{i,a}^*} \abs{f - f_{U_{i,a}^*}}^{\frac{2n}{n-2}} d\mu_\rho \right)^{\frac{n-2}{n}} 
\leq S_c \int_{U_{i,a}^\sharp} \abs{df}^2 dvol,
$$
with $S_c = S_c(n,\kappa)$.
\end{lem}

\begin{preuve}
Set $q =\frac{2n}{n-2}$. For $f \in C^\infty(U_{i,a}^\sharp)$, $i \in \Integer$ :
\begin{eqnarray*}
\int_{U_{i,a}} \abs{f - f_{U_{i,a},\mu_\rho}}^q d\mu_\rho 
&\leq& 2^q \inf_{c \in \Real} \int_{U_{i,a}} \abs{f - c}^q d\mu_\rho \\
&\leq& 2^q \int_{U_{i,a}} \abs{f - f_{U_{i,a},\vol}}^q d\mu_\rho,
\end{eqnarray*} 
so that (\ref{soboann}) (with some small $\delta$ : $0 < \delta < 1 - \kappa^{-1}$) and (\ref{controle}) imply
\begin{eqnarray*}
\left( \int_{U_{i,a}} \abs{f - f_{U_{i,a}}}^q d\mu_\rho \right)^{2/q}
&\leq& \rho(R_{i-1})^{-2/n} C(n,\kappa) \rho(R_{i+1})^{2/n} \int_{U_{i,a}^*} \abs{df}^2 dvol\\
&\leq& C(n,\kappa) \kappa^{2n} \int_{U_{i,a}^*} \abs{df}^2 dvol\\
&\leq& C(n,\kappa)  \int_{U_{i,a}^*} \abs{df}^2 dvol.
\end{eqnarray*}
And such estimates with the pairs $(U_{i,a}^*,U_{i,a}^\sharp)$ also hold for the same reason. 
\end{preuve}


\subsubsection{The discrete Sobolev inequality.}

We consider the weighted graph $(\V,\E,m_\rho)$ associated to the good covering $\U$ of $M$ in $M$, with respect to
$(\mu_\rho,vol)$ (to simplify the notation, we write $m_\rho$ instead of $m_{\mu_\rho}$). What about the structure of the graph ? 
Proposition \ref{litam}, plus the fact that the geometry near $o$ is quasi-euclidian, implies the associated graph, outside a finite subset, 
consists of two trunks, corresponding to neighbourhoods of $o$ and of infinity ; moreover, thanks to the bound $h(n,\kappa)$ on the $h_i$, 
the degrees of the vertices admit an upper bound in terms of $n$ and $\kappa$.    

The measure $m_\rho$ is defined as follows : for each $i \in \Integer$ and $a \in [0,h_i]$,
$$
m_\rho(i,a) = \int_{U_{i,a}} \rho(r)^{-\frac{2}{n-2}} dvol,
$$
so that we can estimate :
$$
\vol(U_{i,a}) \rho(R_{i+1})^{-\frac{2}{n-2}} \leq m_\rho(i,a) \leq  \vol(U_{i,a}) \rho(R_{i-1})^{-\frac{2}{n-2}} ;
$$
using (\ref{bishopgromovdecale}) and (\ref{controle}), this yields
\begin{equation}\label{encadrement}
C(n,\kappa)^{-1} V(R_i) \rho(R_i)^{-\frac{2}{n-2}} \leq m_\rho(i,a) \leq  C(n,\kappa) V(R_i) \rho(R_i)^{-\frac{2}{n-2}}.
\end{equation}
In particular, again with (\ref{bishopgromovdecale}) and (\ref{controle}), this allows us to apply proposition \ref{L1L2} : we are left 
to show that an isoperimetric inequality (\ref{isoperi}) actually occurs. 

Let $\Omega$ be a finite subset of $\V$. Set $l :=\max \set{i \in \Integer, \, \Exists{a}{[0,h_i]}  (i,a) \in \Omega}$. First, we choose a
convenient edge in $\partial \Omega$.

\begin{itemize}
\item If $(l,0)$ belongs to $\Omega$, the edge $e := ((l,0),(l+1,0))$ is in $\partial \Omega$. 
\item Otherwise, we choose $(l,b)$ in $\Omega$. Our choice of $\kappa$ ensures there is a sequence of edges staying on the levels $l$ 
and $l-1$ and which connects $(l,b)$ to $(l,0)$. Among these, there is necessarily an edge which connects a vertex in $\Omega$ to a 
vertex outside $\Omega$ and we call it $e$ : it belongs to $\partial \Omega$. 
\end{itemize}

Then we can write
$$\frac{m_\rho(\Omega)}{m_\rho(\partial \Omega)} \leq \frac{\sum_{i=-\infty}^l \sum_{a=0}^{h_i} m_\rho(i,a)}{m_\rho(e)} 
\leq C(n,\kappa) \sum_{i=-\infty}^l \frac{\sum_{a=0}^{h_i} m_\rho(i,a)}{m_\rho(l,0)}.
$$

With (\ref{controle}), we find  
\begin{eqnarray*}
\frac{m_\rho(\Omega)}{m_\rho(\partial \Omega)} 
&\leq& C(n,\kappa) \sum_{i=-\infty}^l  \frac{V(R_i) \rho(R_i)^{-\frac{2}{n-2}}}{ V(R_l) \rho(R_l)^{-\frac{2}{n-2}}}\\
&\leq& C(n,\kappa) \sum_{i=-\infty}^l  \left[ \frac{V(R_i)}{V(R_l)} \left( \frac{R_i}{R_l} \right)^{-2} \right]^{\frac{n}{n-2}}
\end{eqnarray*}
so that (\ref{volboules1}) gives
\begin{eqnarray*}
\frac{m_\rho(\Omega)}{m_\rho(\partial \Omega)} 
&\leq& C(n,\kappa) C_o^{-\frac{n}{n-2}} \sum_{i=-\infty}^l   \left( \frac{R_i}{R_l} \right)^{\frac{n(\nu-2)}{n-2}} \\
&=& C(n,\kappa) C_o^{-\frac{n}{n-2}} \sum_{j=0}^\infty   \kappa^{- j \frac{n(\nu-2)}{n-2}}\\
&=& \frac{ C(n,\kappa) C_o^{-\frac{n}{n-2} } }{ 1-\kappa^{- \frac{ n(\nu-2) }{ n-2 } }} ,
\end{eqnarray*}
since $\nu>2$.

Then (\ref{L1L2}) and (\ref{isoperisobo}), with $k=\infty$, lead to the
\begin{lem}\label{sobodisc}
For any $1\leq p<\infty$, there exists a constant $S_d$, depending on $p$, $\kappa$, 
$n$, $C_o$, $\nu$, 
such that for every real function $f$ with finite support in $\V$ :
$$ 
\left(\sum_{v \in \V} \abs{f(v)}^{p} m_\rho(v)\right)^{\frac{1}{p}} 
\leq S_d \left(\sum_{(v,w) \in \E} \abs{f(v)-f(w)}^p m_\rho(v,w)\right)^{\frac{1}{p}}.
$$
\end{lem}


\subsubsection{Conclusion.}

\begin{thm}[Weighted Sobolev inequality]\label{inegalite}
Let $M^n$, $n \geq 3$, be a connected complete riemannian manifold with nonnegative Ricci curvature. Assume that there
exists $o \in M$, $\nu > 2$ and  $C_o >0$ such that 
$$
\forall R_2 \geq R_1 > 0, \  \frac{V(o,R_2)}{V(o,R_1)} \geq C_o \left( \frac{R_2}{R_1} \right)^\nu.
$$
Then $M$ satisfies the weighted Sobolev inequality 
$$
\Forall{f}{C^\infty_c(M)} \left( \int_M  \abs{f}^{\frac{2 n}{n-2}} \rho_o(r_o)^{-\frac{2}{n-2}} dvol \right)^{1- \frac{2}{n}} 
\leq S \int_M \abs{df}^2 dvol. 
$$
Here, $S$ can be chosen to depend only on $n$, $C_o$, $\nu$. 
\end{thm}

\begin{preuve}
We just use \ref{recollsobodir}, \ref{sobocont} and \ref{sobodisc}. 
\end{preuve}

\begin{rem}
If one prefers polynomial weights, note \ref{inegalite} implies there is a constant $\tilde{S}$ such that
$$
\Forall{f}{C^\infty_c(M)} \left( \int_M  \abs{f}^{\frac{2 n}{n-2}} \underline{r_o}^{-\frac{2(n-\nu)}{n-2}} dvol \right)^{1- \frac{2}{n}} 
\leq \tilde{S} \int_M \abs{df}^2 dvol, 
$$
where $\underline{r_o}$ is the function which is equal to $1$ inside $B(o,1)$ and to $r_o$ outside this ball (just use (\ref{volboules2})).
Observe we cannot write $r_o$ instead of $\underline{r_o}$, unless $\nu = n$. The obstruction to do this for the Sobolev inequality is 
that locally the weight would not fit : the corresponding inequality is false on $\Real^n$, hence on any riemannian manifold (use the 
family of functions $\max(1-r_o/\epsilon,0)$, $\epsilon>0$). Note also that $\tilde{S}$ depends on  $n$, $C_o$, $\nu$ and $V(o,1)$. 
\end{rem}

Let us introduce some notation for the best constant in our inequality.
\begin{defn}\label{bestsobolev}
Let $M^n$ be a connected complete riemannian manifold, $n \geq 3$. For every $o \in M$, we define the riemannian invariant
$$
S_{o}(M) := \sup_{f \in C^\infty_c(M) \backslash \set{0}} 
 \frac{\left(\int_M  \abs{f}^{\frac{2 n}{n-2}} \rho_o(r_o)^{-\frac{2}{n-2}} dvol \right)^{1- \frac{2}{n}}}{\int_M \abs{df}^2 dvol}.
$$
\end{defn}

The same method gives the

\begin{thm}
Let $M^n$, $n \geq 3$, be a connected noncompact complete riemannian manifold with nonnegative Ricci curvature. Assume that there exists 
$o \in M$, $\nu > 1$ and  $C_o >0$ such that 
$$
\forall R_2 \geq R_1 >0, \  \frac{V(o,R_2)}{V(o,R_1)} \geq C_o \left( \frac{R_2}{R_1} \right)^\nu.
$$
Then if $\beta > -\frac{\nu-2}{n-\nu}$, $M$ satisfies the weighted Sobolev inequality 
$$
\Forall{f}{C^\infty_c(M)} \left( \int_M  \abs{f}^{\frac{2 n}{n-2}} \rho_o(r_o)^{\frac{n\beta - 2}{n-2}}  dvol \right)^{\frac{n-2}{n}} 
\leq S_\beta \int_M \abs{df}^2 \rho_o(r_o)^\beta dvol, 
$$
with $S_\beta=S_\beta(n,C_o,\nu,\beta)$.
\end{thm}

\begin{preuve}
We wish to apply (\ref{recollsobodir}) to the measures $\rho_o(r_o)^{\frac{n\beta - 2}{n-2}} dvol$ and 
$\rho_o(r_o)^\beta dvol$ and the same good
covering. Our choice of weights ensures the continuous Sobolev inequality, as in (\ref{sobocont}) : 
for $i$ in $\Integer$, $a \in [0,h_i]$ and $f$ in $C^\infty(U_i^*)$, (\ref{soboann}) yields 
\begin{eqnarray*}
&\left( \int_{U_{i,a}} \abs{f - f_{U_{i,a}}}^\frac{2n}{n-2} \rho_o(r_o)^{\frac{n\beta - 2}{n-2}} dvol \right)^{1-2/n}& \\
&\leq C(n,\kappa) \rho_o(R_i)^{\frac{n\beta - 2}{n-2}} \rho_o(R_i)^{\frac{2}{n}} \int_{U_{i,a}^*} \abs{df}^2 dvol& \\
&\leq C(n,\kappa) \rho_o(R_i)^{\frac{n\beta - 2}{n}} \rho_o(R_i)^{\frac{2}{n}} \rho_o(R_i)^\beta
\int_{U_{i,a}^*} \abs{df}^2 \rho_o(r_o)^\beta dvol& \\
&= C(n,\kappa) \int_{U_{i,a}^*} \abs{df}^2 \rho_o(r_o)^\beta dvol.&
\end{eqnarray*}

As for the discrete inequality, we proceed as in the proof of \ref{sobodisc}. Essentially, using the same notations as in this proof, we
obtain 
\begin{eqnarray*}
\frac{m(\Omega)}{m(\partial \Omega)} 
&\leq& C(n,\kappa) \sum_{i=-\infty}^l  \frac{V(R_i) \rho(R_i)^{\frac{n\beta - 2}{n-2}}}{ V(R_l) \rho(R_l)^{\frac{n\beta - 2}{n-2}}}\\
&\leq& C(n,\kappa) \sum_{i=-\infty}^l  \left[ \left( \frac{V(R_i)}{V(R_l)} \right)^{1-\beta} 
\left( \frac{R_i}{R_l} \right)^{n\beta - 2} \right]^{\frac{n}{n-2}}
\end{eqnarray*}
so that (\ref{volboules1}) gives
$$\frac{m(\Omega)}{m(\partial \Omega)} 
\leq C(n,\kappa) C_o^{-\frac{n(1-\beta)}{n-2}} \sum_{j=0}^\infty   \kappa^{- j \frac{n(\nu-2+\beta(n-\nu))}{n-2}}
$$
which is finite thanks to our assumption on $\beta$.
\end{preuve}

\begin{rem}
In particular, for $\beta = 1$, the inequality reads
$$
\Forall{f}{C^\infty_c(M)} \left( \int_M  \abs{f}^{\frac{2 n}{n-2}} \frac{r_o^n}{V(o,r)} dvol \right)^{\frac{n-2}{n}} 
\leq S \int_M \abs{df}^2 \frac{r_o^n}{V(o,r)} dvol
$$
The picture is the following : the volume growth of balls is in general not euclidian 
(i.e. it does not behave like $r^n$) and therefore we cannot hope
to find a nonweighted Sobolev inequality (cf. next paragraph) ; nevertheless, by radially modifying the riemannian measure so that 
it has euclidian growth, we manage to obtain a Sobolev inequality.
\end{rem}


\subsubsection{What does a weighted Sobolev inequality implies on the volume growth of balls ?}

\begin{prop}
Let $M^n$, $n \geq 3$, be a connected noncompact complete riemannian manifold with nonnegative Ricci curvature. Assume that there
exists $o \in M$, $\alpha \geq 0$ and  $S >0$ such that 
$$
\Forall{f}{C^\infty_c(M)} \left( \int_M  \abs{f}^{\frac{2 n}{n-2}} \underline{r_o}^{-\alpha} dvol \right)^{\frac{n-2}{n}} 
\leq S \int_M \abs{df}^2 dvol.
$$
Then there is a constant $A_o >0$ such that  
$$
\forall R\geq 1, \  V(o,R) \geq A_o R^\nu,
$$
where $\nu$ is the real number defined by $\alpha = 2 \frac{n-\nu}{n-2}$.
\end{prop}

\begin{preuve}
As usual, we set $q = 2n/(n-2) > 2$. 
Then we fix $R \geq 2$ and $0 < t \leq R/2$ and consider the Lipschitz function
$$
f := \max(t - d(.,S(o,R)), 0)) :
$$
$f=t$ on the sphere $S(o,R)$, $f=0$ outside some $t$-neighbourhood of this sphere and, on this $t$-neighborhood, 
it decreases radially at unit speed. Thus
$$
\int_M  \abs{f}^{q} r^{-\alpha} dvol \geq (t/2)^q (R+t)^{-\alpha} \vol(A(R-t/2,R+t/2)
$$
and
$$
\int_M \abs{df}^2 dvol \leq \vol(A(R-t,R+t).
$$
The Sobolev inequality yields :
$$
(t/2)^2 (R+t)^{-2\alpha/q} \vol(A(R-t/2,R+t/2)^{2/q} \leq S \vol(A(R-t,R+t).
$$
For $i \in \Natural^*$, we apply this to $t=2^{-i} R$. With $V_i := \vol(A(R(1-2^{-i}),R(1+2^{-i}))$, 
$$
R^2 4^{-i-1} ((1+2^{-i})R)^{-2\alpha/q} V_{i+1}^{2/q} \leq S V_i.
$$
By induction, there is a constant $C$ which does not depend on $R$ such that for every $i \geq 1$
$$
\vol(B(2R)) \geq V_1
\geq   \left( C R^{2-2\alpha/q} \right)^{\sum_{j=0}^{i-1} (2/q)^j} 
\left( \prod_{j=0}^{i-1}  (4^{-j})^{(2/q)^j} \right) V_i.
$$
As a riemannian manifold is locally quasi-euclidian, for $i \To \infty$, 
$$
V_i^{(2/q)^i} \geq  \left( C(R) (2^{-i}R)^n \right)^{(2/q)^i}  \To 1.
$$
Eventually, 
$$
\vol(B(2R)) \geq C^{\frac{1}{1-2/q}} R^{\frac{2-2\alpha/q}{1-2/q}} \prod_{j=0}^{\infty}  (4^{-j})^{(2/q)^j}.
$$
And indeed, $\nu = \frac{2-2\alpha/q}{1-2/q}$ is the same as
$\alpha =  2 \frac{n-\nu}{n-2}$.
\end{preuve}


\subsection{The Hardy inequality.}

With \ref{recollsobodir}, we can also patch local Poincaré inequalities together. Working under the same assumptions as above, the global inequality we find is a Hardy inequality. 

\begin{thm}[Hardy inequality]\label{hardy}
Let $M^n$, $n \geq 3$, be a connected noncompact complete riemannian manifold with nonnegative Ricci curvature. Fix some $p\geq 1$. Assume 
that there exists $o \in M$, $\nu > p$ and  $C_o >0$ such that 
$$
\forall R_2 > R_1 \geq 1, \  \frac{V(o,R_2)}{V(o,R_1)} \geq C_o \left( \frac{R_2}{R_1} \right)^\nu.
$$
Then $M$ satisfies the Hardy inequality 
$$
\Forall{f}{C^\infty_c(M)} \int_M  \abs{f}^{p} r_o^{-p}  dvol
\leq H \int_M \abs{df}^p dvol, 
$$
with a constant $H$ depending only on $n$, $C_o$, $\nu$, $p$.
\end{thm}

\begin{preuve}
The proof consists in applying \ref{recollsobodir} with $k=\infty$. We will use the same "good" covering $\U$ as in paragraph \ref{goodcov}, 
noticing it is also "good" for the pair of measures $(r^{-p} dvol,dvol)$. 

We need a continuous Poincaré inequality. Indeed, as for \ref{sobocont}, if we choose $i \in \Integer$ and $a \in [0,h_i]$, each smooth 
function $f$ on $U_{i,a}^\sharp$ satisfies
\begin{eqnarray*}
\int_{U_{i,a}} \abs{f - f_{U_{i,a},\mu_\alpha}}^p r^{-p} dvol 
&=& \inf_{c \in \Real} \int_{U_{i,a}} \abs{f - c}^p r^{-p} dvol \\
&\leq& \int_{U_{i,a}} \abs{f - f_{U_{i,a}}}^p r^{-p} dvol
\end{eqnarray*} 
so that, with \ref{poincann},  
\begin{eqnarray*}
\int_{U_{i,a}} \abs{f - f_{U_{i,a}}}^p r^{-p} dvol 
&\leq& C(n,\kappa) R_{i-1}^{-p} R_{i+1}^p \int_{U_{i,a}} \abs{f - f_{U_{i,a}}}^p r^{-p} dvol\\
&\leq& C(n,\kappa) \int_{U_{i,a}^*} \abs{df}^p dvol.
\end{eqnarray*} 
And the same argument works with the pairs $(U_{i,a}^*,U_{i,a}^\sharp)$.

The discrete inequality required in \ref{recollsobodir} follows from the argument of \ref{sobodisc} ; here, we estimate the discrete 
isoperimetric quotient by
$$
C(n,\kappa) \sum_{i=-\infty}^l   \frac{V(R_i)}{V(R_l)} \left( \frac{R_i}{R_l} \right)^{-p} 
$$
which is bounded by
$$
C(n,\kappa) C_o \sum_{j=0}^\infty   \kappa^{- j (\nu-p)} <\infty
$$
thanks to our assumption on the volume growth of balls.
\end{preuve}

\vspace{\baselineskip}

For convenience, we give a name to the best constant in the Hardy inequalities.

\begin{defn}\label{besthardy}
Let $M$ be a connected complete riemannian manifold. For $o \in M$ and $r_o:=d(o,.)$, we define the riemannian invariant
$$
H_{o}(M) := \sup_{f \in C^\infty_c(M) \backslash \set{0}} 
\frac{\int_M  \abs{f} r_o^{-1} dvol}{\int_M \abs{df} dvol}.
$$
\end{defn}


\section{Weighted Sobolev inequalities and Schrödinger operators.}

In this section, we explain a few analytical consequences of the weighted Sobolev inequality. They will find geometric applications in the
last section. We assume here that $M^n$ is a connected complete noncompact manifold such that for some point $o$ in $M$ and $S > 0$, 
the following weighted Sobolev inequality is true :
$$
\Forall{f}{C^\infty_c(M)} \left( \int_M \abs{f}^{\frac{2 n}{n-2}} \rho_o(r_o)^{-\frac{2}{n-2}} dvol \right)^{1-\frac{2}{n}} 
\leq S \int_M \abs{df}^2 dvol.
$$ 
As previously, we will often write $\rho(r)$ for $\rho_o(r_o)$, but also 
$$
d\mu_\rho = \rho(r)^{-\frac{2}{n-2}} dvol
$$
and
$$
q = \frac{2n}{n-2}.
$$

We consider a smooth euclidian vector bundle $E \To M$, endowed with a compatible connection $\nabla$. 
We will always denote by $(.)$ the pointwise scalar product on a euclidian vector bundle, by $\abs{.}$ the pointwise norm, 
by $\overline{\Delta} = \nabla^* \nabla$ the Bochner laplacian (or "rough laplacian"). 
Our interest lies in Schrödinger operators $\overline{\Delta} + V$, where $V$ is a continuous field of symmetric endomorphisms of $E$. 
We decompose $V$ as $V=V_+ - V_-$, where $V_+$ and $V_-$ are fields of positive symmetric endomorphisms of $E$. 
We describe here some consequences of the weighted Sobolev inequality on these operators. 

\vspace{\baselineskip}


\subsection{A vanishing theorem.}

The following theorem is a generalization of \cite{Car1}.

\begin{thm}[Vanishing theorem]\label{annulation}
Fix $m>1$ and assume the potential $V$ satisfies
$$
S \,\, \left(\int_M \abs{V_-}^{\frac{n}{2}} \rho(r) dvol \right)^{\frac{2}{n}} < \epsilon(m),
$$
where 
$$
\epsilon(m) = \left\{ \begin{array}{ll}
\frac{2}{m} & \text{if $m\geq2$,} \\
\frac{2}{m} \left( 2 - \frac{2}{m} \right) & \text{if $1<m<2$,}
\end{array}\right.
$$
Then every locally Lipschitz section $\sigma$ of $E$ such that
$$
\int_{A(R/2,R)} \abs{\sigma}^m dvol = o(R^{2})
$$
and
$$
(\overline{\Delta} \sigma + V \sigma,\sigma) \leq 0
$$
is identically zero.
\end{thm}


\begin{rem}
In this statement, the distribution $(\overline{\Delta} \sigma,\sigma)$ is defined by : \\
$\displaystyle{\Forall{\phi}{C^\infty_c(M)} <(\overline{\Delta} \sigma,\sigma),\phi> = \int_M (\nabla \sigma,\nabla(\phi \sigma)) dvol}$.
\end{rem}

\begin{demo}
We first treat the case $m = 2$. Let $R$ be a positive number. Let us choose a smooth function $\chi$ which is $1$ on $B(R)$, 
$0$ on $M \backslash B(2R)$, takes its values in $[0,1]$ and satisfies $\abs{d \chi} \leq 2/R$.
We apply the weighted Sobolev inequality to the Lipschitz function $\chi \abs{\sigma}$ (we omit the riemannian measure in the next
formulas so as to make them easier to read) :
$$
\left( \int_M \chi^q \abs{\sigma}^{q} \rho(r)^{-\frac{2}{n-2}} \right)^{\frac{2}{q}}
\leq S \int_M \abs{d(\abs{\chi \sigma})}^2
\leq S \int_M \abs{\nabla (\chi \sigma)}^2,
$$
where we used the Kato inequality.
Now, 
$$
\overline{\Delta} (\chi \sigma) = \chi \overline{\Delta} \sigma + (\Delta \chi ) \sigma - 2 \nabla_{\grad \chi } \sigma
$$ 
and integration by parts gives 
\begin{eqnarray*}
\int_M \abs{\nabla(\chi \sigma)}^2 &=& \int_M (\overline{\Delta}(\chi \sigma), \chi \sigma)  \\
&=& \int_M \chi^2 (\sigma, \overline{\Delta} \sigma)  + \int_M \chi \Delta \chi \abs{\sigma}^2  
- \frac{1}{2} \int_M (d(\chi^2),d(\abs{\sigma}^2))  \\
&=& \int_M \chi^2 (\sigma, \overline{\Delta} \sigma)  + \int_M \abs{d\chi}^2 \abs{\sigma}^2 \\
&\leq& -\int_M \chi^2 (V_- \sigma,\sigma)  + \frac{4}{R^2} \int_{A(R,2R)} \abs{\sigma}^2 .
\end{eqnarray*}
The Hölder inequality implies
\begin{eqnarray*}
-\int_M \chi^2 (V_- \sigma,\sigma)  &\leq& \int_M \chi^{2} \abs{\sigma} \rho(r)^{-\frac{2}{n}} 
\abs{V_-} \rho(r)^{\frac{2}{n}}   \\
&\leq& \left(\int_M \chi^q \abs{\sigma}^{q} \rho(r)^{-\frac{2}{n-2}} \right)^{\frac{2}{q}}
 \underbrace{\left(\int_M \abs{V_-}^{\frac{n}{2}} \rho(r)  \right)^{\frac{2}{n}}}_{N_V}.
\end{eqnarray*}
All in all, we find
$$
(1 - S N_V) \left( \int_M \chi^q \abs{\sigma}^{q} \rho(r)^{-\frac{2}{n-2}}  \right)^{\frac{2}{q}}
\leq \frac{4 S}{R^2} \int_{A(R,2R)} \abs{\sigma}^{2} 
$$
and the assumption on the potential allows us to write
$$
\left( \int_{B(R)} \abs{\sigma}^{q} \rho(r)^{-\frac{2}{n-2}}  \right)^{\frac{2}{q}}
\leq  \frac{4 S}{1 - S N_V} \frac{1}{R^2} \int_{A(R,2R)} \abs{\sigma}^{2}.
$$
Letting $R \To +\infty$, we obtain $\sigma=0$.

Now, we turn to the case $m \geq 2$. First note that the Kato inequality implies
$$
\abs{\sigma} \Delta \abs{\sigma} = \abs{d \abs{\sigma}}^2 + \frac{1}{2} \Delta \abs{\sigma}^{2} 
\leq \abs{\nabla \sigma}^{2} + \frac{1}{2} \Delta \abs{\sigma}^{2},
$$
and since
$$
(\sigma, \overline{\Delta} \sigma) = \abs{\nabla \sigma}^2 + \frac{1}{2} \Delta \abs{\sigma}^2,
$$ 
this means we always have the inequality
$$
\abs{\sigma} \Delta \abs{\sigma} \leq (\sigma, \overline{\Delta} \sigma)
$$ 
and in our setting
$$
\abs{\sigma}\Delta \abs{\sigma} \leq (\sigma, V_- \sigma) \leq \abs{V_-} \abs{\sigma}^{2}.
$$
So, if $u:=\abs{\sigma}^{m/2}$, 
\begin{eqnarray*}
u\Delta u &=& \abs{\sigma}^{m/2} \Delta \abs{\sigma}^{m/2} \\
&=& \frac{m}{2} \abs{\sigma}^{m-2} \abs{\sigma} \Delta \abs{\sigma} 
- \frac{m}{2} \left( \frac{m}{2} -1 \right) \abs{\sigma}^{m-2} \abs{d\abs{\sigma}}^2\\
&\leq& \frac{m}{2} \abs{\sigma}^{m} \abs{V_-} 
- \frac{m}{2} \left( \frac{m}{2} -1 \right) \abs{\sigma}^{m-2} \abs{d\abs{\sigma}}^2\\
&=& \frac{m}{2} u^2 \abs{V_-}
- \left( 1- \frac{2}{m} \right) \abs{d\abs{u}}^2.
\end{eqnarray*}
Thus $u(\Delta u - \frac{m}{2} \abs{V_-} u) \leq 0$ and we can apply the case $m=2$ to $u$.

For the case $1<m<2$, we go back to the weighted Sobolev inequality, applied to the locally Lipschitz function $\chi u_\epsilon^{m/2}$, 
where $u_\epsilon = \sqrt{\abs{\sigma}^2 + \epsilon}$, $\epsilon>0$ : 
\begin{eqnarray*}
&\frac{1}{S} \left( \int_M \chi^q u_\epsilon^{\frac{mq}{2}} \rho(r)^{-\frac{2}{n-2}} \right)^{\frac{2}{q}}&\\
&= \int_M \abs{d\chi}^2 u_\epsilon^{m}  + \int_M \chi^2 \abs{ d (u_\epsilon^{m/2})}^2  + 2 \int_M (u_\epsilon^{m/2} d\chi, \chi
d(u_\epsilon^{m/2}))& \\
&\leq (1+1/b) \int_M \abs{d\chi}^2 u_\epsilon^{m}  +  (1+b) \int_M \chi^2 \abs{ d (u_\epsilon^{m/2})}^2&
\end{eqnarray*}
for any $b>0$.
Integration by parts yields
\begin{eqnarray*}
\int_M &\chi^2& \abs{d(u_\epsilon^{m/2})}^2 \\
= &\int_M& (\chi^2 d(u_\epsilon^{m/2}),d(u_\epsilon^{m/2})) \\ 
= &\int_M& 2\chi (u_\epsilon^{m/2} d\chi, d(u_\epsilon^{m/2})) + \int_M \chi^2 u_\epsilon^{m/2} \Delta (u_\epsilon^{m/2}) \\ 
= &2& \int_M (u_\epsilon^{m/2} d\chi, \chi d(u_\epsilon^{m/2})) + \frac{m}{2} \int_M \chi^2 u_\epsilon^{m-1} \Delta u_\epsilon \\
  &+& \left( \frac{2}{m} - 1 \right) \int_M \chi^2 \abs{d(u_\epsilon^{m/2})}^2.
\end{eqnarray*}
So, if $a>0$, 
\begin{eqnarray*}
\int_M \chi^2 \abs{d(u_\epsilon^{m/2})}^2  
\leq (\frac{2}{m} - 1 + a) &\int_M& \chi^2 \abs{d(u_\epsilon^{m/2})}^2\\
 + \frac{m}{2} &\int_M& \chi^2 u_\epsilon^{m-1} \Delta u_\epsilon + \frac{1}{a} \int_M \abs{d\chi}^2 u_\epsilon^{m}
\end{eqnarray*}
and if moreover $a<2-2/m$,
$$
\int_M \chi^2 \abs{d(u_\epsilon^{m/2})}^2  \leq 
(2 - \frac{2}{m} - a)^{-1} \left(  \frac{m}{2} \int_M \chi^2 u_\epsilon^{m-1} \Delta u_\epsilon 
+ \frac{1}{a} \int_M \abs{d\chi}^2 u_\epsilon^{m} \right).
$$
Thus
$$
\frac{1}{S} \left( \int_M \chi^q u_\epsilon^{\frac{mq}{2}} \rho(r)^{-\frac{2}{n-2}} \right)^{\frac{2}{q}}
\leq C(m,a,b) \int_M \abs{d\chi}^2 u_\epsilon^{m} +
D(m,a,b) \int_M \chi^2 u_\epsilon^{m-1} \Delta u_\epsilon  
$$
where 
$$
C(m,a,b) = 1+ 1/b + \frac{1+b}{a(2-2/m -a)}$$ 
and 
$$
D(m,a,b)= \frac{(1+b)m}{2(2-2/m -a)}.$$
We compute
$$
u_\epsilon \Delta u_\epsilon = (\sigma, \overline{\Delta} \sigma) -  \frac{\epsilon\abs{\nabla \sigma}^2}{u_\epsilon^2}
- \frac{\abs{\sigma}^2 \abs{\nabla \sigma}^2 - (\sigma,\nabla \sigma)^2}{u_\epsilon^2}, 
$$
to ensure
$$
u_\epsilon \Delta u_\epsilon \leq (\sigma, \overline{\Delta} \sigma) \leq \abs{V_-} \abs{\sigma}^2.
$$
Therefore,
$$
\frac{1}{S} \left( \int_M \chi^q u_\epsilon^{\frac{mq}{2}} \rho(r)^{-\frac{2}{n-2}} \right)^{\frac{2}{q}}
\leq C(m,a,b) \int_M \abs{d\chi}^2 u_\epsilon^{m} +
D(m,a,b) \int_M \chi^2 u_\epsilon^{m-2} \abs{V_-} \abs{\sigma}^2 
$$
and when $\epsilon \To 0$,
$$
\frac{1}{S} \left( \int_M \chi^q \abs{\sigma}^{\frac{mq}{2}} \rho(r)^{-\frac{2}{n-2}} \right)^{\frac{2}{q}}
\leq C(m,a,b) \int_M \abs{d\chi}^2 \abs{\sigma}^{m} +
D(m,a,b) \int_M \chi^2 \abs{V_-} \abs{\sigma}^{m}.   
$$
As previously, this implies :
$$
\left( \int_{B(R)} \abs{\sigma}^{mq/2} \rho(r)^{-\frac{2}{n-2}}  \right)^{\frac{2}{q}}
\leq \frac{1}{1 - S N_V D(m,a,b)} \frac{4 S C(m,a,b)}{ R^2} \int_{A(R,2R)} \abs{\sigma}^{m},
$$
providing
$$
N_V < \frac{1}{S D(m,a,b)} = \frac{2}{m S} \frac{1}{1 + b} (2 - 2/m -a),
$$
which, under our assumption on $V$, can always be achieved by choosing sufficiently small $a$ and $b$. Letting $R \To \infty$, we prove the
claim.
\end{demo}

\vspace{\baselineskip}


\subsection{Some general decay estimates.}\label{decaylemmas}

Now what can we say if we only have
$$
\int_M \abs{V_-}^{\frac{n}{2}} \rho(r) dvol < \infty \,?
$$ 
Adapting a technique developped in \cite{BKN}, we can prove some decay estimates on the sections $\sigma$ such that 
$\overline{\Delta} \sigma + V \sigma \leq 0$. We prove three general lemmas and we will see later (\ref{ricciflat}) how to apply them 
in a geometric setting, where the potential depends on the section $\sigma$. The idea is to implement a Nash-Moser iteration : 
this is the third lemma. But this lemma only works under a technical assumption on the potential, which can be ensured by the first lemma.
Finally, the second lemma is a key to a "self-improvement" of the decay estimate we will find. 

\begin{lem}[Initiation]\label{lemdec1}
We assume the potential $V$ satisfies 
$$
\int_M \abs{V_-}^{n/2} \rho(r) dvol < + \infty.
$$
and we consider a locally Lipschitz section $\sigma$ of $E$ such that for some $m>1$ 
$$
\int_{A(R,2R)} \abs{\sigma}^m dvol = o(R^2)
$$
and
$$
(\sigma,\overline{\Delta} \sigma + V \sigma) \leq 0.
$$
Then for large $R$ :
$$
\left(\int_{M \backslash B(2R)} \abs{\sigma}^{\frac{m q}{2}} d\mu_\rho \right)^{\frac{2}{q}}
\leq \frac{C}{R^{2}} \int_{A(R,2R)} \abs{\sigma}^m dvol.
$$
\end{lem}

\begin{preuve}
Proceeding as in the proof of the vanishing theorem, we find for $u:=\abs{\sigma}^{m/2}$ and $\chi \in C^\infty_c(M)$ :
$$
\left( \int_M \chi^q u^{q} \rho(r)^{-\frac{2}{n-2}} \right)^{\frac{2}{q}}
\leq C \left( \int_M \chi^2 u^2 \abs{V_-}  + \int_M \abs{d\chi}^2 u^2  \right),
$$
and, using Hölder inequality, this yields :
\begin{eqnarray*}
\left( \int_M \chi^q u^{q} \rho(r)^{-\frac{2}{n-2}} \right)^{\frac{2}{q}}
\leq &C& \left( \int_{\supp \chi} \abs{V_-}^{\frac{n}{2}} \rho(r)  \right)^{\frac{2}{n}}   
\left( \int_M \chi^q u^{q} \rho(r)^{-\frac{2}{n-2}} \right)^{\frac{2}{q}} \\
+ &C& \int_M \abs{d\chi}^2 u^2 .
\end{eqnarray*}

Now we set $R>>1$, $R'> 2R$ and we choose $\chi$ with support in $A(R,2R')$, with value $1$ on $[2R,R']$, satisfying
$\abs{d\chi} \leq \frac{2}{R}$ on $A(R,2R)$ and $\abs{d\chi} \leq \frac{2}{R'}$ on $A(R',2R')$. Thus
\begin{eqnarray*}
\left( \int_{M} \chi^q u^{q} \rho(r)^{-\frac{2}{n-2}} \right)^{\frac{2}{q}}
\leq &C& \left( \int_{A(R,2R')} \abs{V_-}^{\frac{n}{2}} \rho(r)  \right)^{\frac{2}{n}}   
\left( \int_{M} \chi^q u^{q} \rho(r)^{-\frac{2}{n-2}} \right)^{\frac{2}{q}} \\
&+& \frac{C}{R^2} \int_{A(R,2R)} u^2  + \frac{C}{R'^2} \int_{A(R',2R')} u^2 .
\end{eqnarray*}

By assumption, the integral $\int_M \abs{V_-}^{n/2} \rho(r) $ is finite : we can make the quantity
$\int_{B(R)^c} \abs{V_-}^{n/2} \rho(r) $ as small as we like, by choosing a large $R$, so that
$$
\left( \int_{A(2R,R')} u^{q} \rho(r)^{-\frac{2}{n-2}} \right)^{\frac{2}{q}}
\leq \frac{C}{R^2} \int_{A(R,2R)} u^2  + \frac{C}{R'^2} \int_{A(R',2R')} u^2 .
$$
Letting $R' \To \infty$ we find
\begin{eqnarray*}
\left(\int_{M \backslash B(2R)} u^{q} \rho(r)^{-\frac{2}{n-2}} \right)^{\frac{2}{q}}
&\leq& \frac{C}{R^2} \int_{A(R,2R)} u^{2} .
\end{eqnarray*}
\end{preuve}

\vspace{\baselineskip}

\begin{lem}[Key to the self-improvement]\label{lemdec2}
We assume the potential $V$ satisfies
$$
\int_M \abs{V_-}^{n/2} \rho(r) dvol < + \infty
$$
and we consider a locally Lipschitz section $\sigma$ of $E$ belonging to $L^m(E,\mu_\rho)$ for some $m>q/2$,
such that
$$
(\sigma,\overline{\Delta} \sigma + V \sigma) \leq 0.
$$
Then for large $R$ :
$$
\int_{M \backslash B(2R)} \abs{\sigma}^{m} d\mu_\rho
\leq C \int_{A(R,2R)} \abs{\sigma}^{m} d\mu_\rho.
$$
As a consequence,
$$
\int_{M \backslash B(R)} \abs{\sigma}^{m} d\mu_\rho =O(R^{-a}),
$$
for some a>0. 
\end{lem}

\begin{rem}
The proof will show that $a$ can be chosen so that it depends continuously on $m$.
\end{rem}

\vspace{\baselineskip}

\begin{preuve}
Set $m' = 2m/q$. The preceding proof says that for large $R$, with the same truncature function $\chi$ and $u:=\abs{\sigma}^{m'/2}$ :
$$
\left( \int_{M} \chi^q u^{q} \rho(r)^{-\frac{2}{n-2}} \right)^{\frac{2}{q}}
\leq  C \int_{M}  \abs{d\chi}^2 u^{2} .
$$
we again use the Hölder inequality :
$$
\left( \int_{A(2R,R')} u^{q} \rho(r)^{-\frac{2}{n-2}} \right)^{\frac{2}{q}}
\leq C \left( \int_{M} \abs{d\chi}^n \rho(r)  \right)^{\frac{2}{n}}
\left( \int_{\supp d\chi} u^{q} \rho(r)^{-\frac{2}{n-2}} \right)^{\frac{2}{q}}.
$$

Now,
$$
\int_{A(R,2R)} \abs{d\chi}^n \rho(r) 
\leq C R^{-n} \rho(2R) \vol A(R,2R)
\leq C 
$$
and also
$$
\int_{A(R',2R')} \abs{d\chi}^n \rho(r)  \leq C,
$$
so that
$$
\left( \int_{A(2R,R')} \abs{\sigma}^{m} \rho(r)^{-\frac{2}{n-2}} \right)^{\frac{2}{q}}
\leq C \left( \int_{A(R,2R) \cup A(R',2R')} \abs{\sigma}^{m} \rho(r)^{-\frac{2}{n-2}} \right)^{\frac{2}{q}}.
$$
Letting $R' \To \infty$, we find the first part of the claim :
$$
\int_{M \backslash B(2R)} \abs{\sigma}^{m} \rho(r)^{-\frac{2}{n-2}} 
\leq C \int_{A(R,2R)} \abs{\sigma}^{m} \rho(r)^{-\frac{2}{n-2}}.
$$

Set $I(R)= \int_{M \backslash B(R)} \abs{\sigma}^{m} \rho(r)^{-\frac{2}{n-2}}$. We proved that for large $R$,
$$
I(2R) \leq C (I(R) - I(2R)),
$$
i.e.
$$
I(2R) \leq \frac{C}{C+1} I(R).
$$
Fix a large $R_1$ and denote by $k_R$ the integer such that
$$
\log_2 R/R_1 \leq k_R < \log_2 2R/R_1. 
$$
Iterating the inequality, we find
$$
I(R) \leq \left( \frac{C}{C+1} \right)^{k_R} I(R/2^{k_R}) \leq  \left( \frac{C}{C+1} \right)^{k_R} \norm{\sigma}_{L^m(E,\mu_\rho)}^m
$$ 
so
$$
I(R) \leq C \left( \frac{C}{C+1} \right)^{\log_2 R} = C R^{\log_2 \left( \frac{C}{C+1} \right)},
$$
hence the second statement, since $\frac{C}{C+1} <1$.
\end{preuve}

\vspace{\baselineskip}

\begin{lem}[Nash-Moser iteration]\label{lemdec3}
We assume the potential $V$ satisfies, for some $x> n/2$, 
$$
\left( \int_{A(R,2R)} \abs{V_-}^x \rho(r)^{\frac{x-1}{n/2-1}} dvol \right)^\frac{1}{x-n/2}
= O\left( \rho(R)^\frac{2}{n-2}  R^{-2} \right)
$$
and we consider a locally Lipschitz section $\sigma$ in $L^m(M,\mu_\rho)$ for some $m>1$, such that
$$
(\sigma,\overline{\Delta} \sigma + V \sigma) \leq 0.
$$
Then there is a constant $C$ such that for large $R$,
$$
\sup_{A(R,2R)} \abs{\sigma} 
\leq C \left(\rho(R)^\frac{2}{n-2} R^{-2} \right)^{\frac{n}{2m}} \left( \int_{A(R/2,5R/2)} \abs{\sigma}^{m} d\mu_\rho \right)^{1/m}.
$$
\end{lem}

\vspace{\baselineskip}

\begin{preuve}
Fix $\beta \geq m$.
Again with the same technique, one sees that for $\chi \in C^\infty_c(M)$,
\begin{equation}\label{debut}
\left( \int_M \chi^q \abs{\sigma}^{\frac{q \beta}{2}} \rho(r)^{-\frac{2}{n-2}} \right)^{\frac{2}{q}}
\leq C \beta \int_M \chi^2 \abs{\sigma}^{\beta} \abs{V_-}  + C \int_M \abs{d\chi}^2 \abs{\sigma}^\beta .
\end{equation}
In this proof, $C$ denotes a constant which does not depend on $\beta$.

The Hölder inequality implies that for real numbers $t$ and $s$ satisfying 
$$
\frac{1}{x} + \frac{1}{s} + \frac{1}{t} = 1
$$
and
\begin{equation}\label{relhold}
\frac{q}{2s} + \frac{1}{t} = 1,
\end{equation}
we have the estimate 
\begin{eqnarray*}
&\beta& \int_M \chi^2 \abs{\sigma}^{\beta} \abs{V_-}  \\
&\leq& \beta \left( \int_{\supp \chi} \abs{V_-}^{x} \rho(r)^{\frac{x-1}{n/2-1}}  
\right)^{\frac{1}{x}} \left(\int_{M} \chi^q \abs{\sigma}^{\frac{q \beta}{2}} \rho(r)^{-\frac{2}{n-2}} \right)^{\frac{1}{s}}
\left(\int_{M} \chi^2 \abs{\sigma}^{\beta} \rho(r)^{-\frac{2}{n-2}} \right)^{\frac{1}{t}}.
\end{eqnarray*}
Note $t= \frac{x}{x-n/2}$.

The Young inequality, with (\ref{relhold}), yields for each $\epsilon>0$ a constant $C_\epsilon$ such that
\begin{eqnarray*}
\beta \int_M \chi^2 \abs{\sigma}^{\beta} \abs{V_-}  
\leq &\epsilon& \left(\int_{M} \chi^q \abs{\sigma}^{\frac{q \beta}{2}} \rho(r)^{-\frac{2}{n-2}} \right)^{\frac{2}{q}}\\ 
&+& C_\epsilon \beta^t  \left( \int_{\supp \chi} \abs{V_-}^{x} \rho(r)^{\frac{x-1}{n/2-1}}  
\right)^{\frac{t}{x}} \left(\int_{M} \chi^2 \abs{\sigma}^{\beta} \rho(r)^{-\frac{2}{n-2}} \right).
\end{eqnarray*}
Consequently, for small $\epsilon$ (regardless of $\beta$), we obtain in (\ref{debut}) :
\begin{eqnarray*}
\left( \int_M \chi^q \abs{\sigma}^{\frac{q \beta}{2}} \rho(r)^{-\frac{2}{n-2}} \right)^{\frac{2}{q}}
&\leq& C \beta^t  \left( \int_{\supp \chi} \abs{V_-}^{x} \rho(r)^{\frac{x-1}{n/2-1}}  
\right)^{\frac{t}{x}} \left(\int_{M} \chi^2 \abs{\sigma}^{\beta} \rho(r)^{-\frac{2}{n-2}} \right)\\
 &+& C \int_M \abs{d\chi}^2 \abs{\sigma}^\beta .
\end{eqnarray*}

Now we consider truncature functions $\chi$ which, given large $R_1 < R_2 < 5 R_1$ and $0 <\delta \leq R_1/2$,
are equal to $1$ on $A(R_1,R_2)$, are equal to $0$ outside $A(R_1-\delta,R_2+\delta)$ and such that there differential is bounded by
$2/\delta$. Note that our assumption, thanks to (\ref{controle}), implies 
$$
\left( \int_{A(R_1-\delta,R_2+\delta)} \abs{V_-}^x \rho(r)^{\frac{x-1}{n/2-1}} \right)^\frac{t}{x}
\leq C \rho(R_1-\delta)^\frac{2}{n-2}  (R_1-\delta)^{-2}
$$
With this in mind, our estimate gives 
\begin{eqnarray*}
& & \left( \int_{A(R_1,R_2)} \abs{\sigma}^{\frac{q \beta}{2}} \rho(r)^{-\frac{2}{n-2}} \right)^{\frac{2}{q}}\\
\leq &C \beta^t&  \left( \int_{A(R_1-\delta,R_2+\delta)} \abs{V_-}^{x} \rho(r)^{\frac{x-1}{n/2-1}} 
\right)^{\frac{t}{x}} \int_{A(R_1-\delta,R_2+\delta)} \abs{\sigma}^{\beta} \rho(r)^{-\frac{2}{n-2}}\\
 &+& C \rho(R_2+\delta)^\frac{2}{n-2}  \delta^{-2} \int_{A(R_1-\delta,R_2+\delta)} \abs{\sigma}^\beta \rho(r)^{-\frac{2}{n-2}}\\
\leq &C \beta^t&  \rho(R_1-\delta)^\frac{2}{n-2}  (R_1-\delta)^{-2} \int_{A(R_1-\delta,R_2+\delta)} \abs{\sigma}^\beta \rho(r)^{-\frac{2}{n-2}} \\
 &+& C \rho(R_2+\delta)^\frac{2}{n-2}  \delta^{-2} \int_{A(R_1-\delta,R_2+\delta)} \abs{\sigma}^\beta \rho(r)^{-\frac{2}{n-2}} \\
\leq &C \beta^t&  \rho(R_2)^\frac{2}{n-2} \delta^{-2}  
\int_{A(R_1-\delta,R_2+\delta)} \abs{\sigma}^\beta \rho(r)^{-\frac{2}{n-2}},
\end{eqnarray*}
so that, with respect to the measure $\mu_\rho$,
\begin{equation}\label{boucle}
\norm{\sigma}_{L^{\beta q/2} (A(R_1,R_2))} \leq \left(C \beta^t  \rho(R_2)^\frac{2}{n-2} \delta^{-2}\right)^{1/\beta} \norm{\sigma}_{L^{\beta} (A(R_1 - \delta ,R_2 + \delta))}.
\end{equation}

Given some large $R>0$, we set for every $k \in \Natural$ :
$$
\beta_k := m \left (\frac{q}{2} \right )^{k}
\text{, }
\delta_k := 2^{-k-1} R,
$$
$$
R_{1,k} := R - \sum_{i=1}^k \delta_i
\text{, }
R_{2,k} :=  2R + \sum_{i=1}^k \delta_i.
$$
Iterating (\ref{boucle}), we find
$$
\norm{\sigma}_{L^{\beta_k} (A(R,2R))} \leq C_k \norm{\sigma}_{L^{\beta_0} (A(R_{1,k},R_{2,k}) )},
$$
where the constant is estimated by 
\begin{eqnarray*}
C_k &\leq& \prod_{i=0}^{k-1} \left(C \beta_i^t  \rho(R)^\frac{2}{n-2} R^{-2} 4^i \right)^{1/\beta_i} \\
&\leq& \left(C \rho(R)^\frac{2}{n-2} R^{-2} \right)^{\sum_{i=0}^{k-1} 1/\beta_i} 
\left(4 (q/2)^t \right)^{\sum_{i=0}^{k-1} i/\beta_i}.
\end{eqnarray*}

Since $\sum_{i=0}^\infty \frac{1}{\beta_i} = \frac{n}{2m}$ et $\sum_{i=0}^\infty \frac{i}{\beta_i} < \infty$,
we find :
$$
\overline{\lim_{k \To \infty}} C_k \leq C \left(\rho(R)^\frac{2}{n-2} R^{-2} \right)^{\frac{n}{2m}},
$$
so
$$
\sup_{A(R,2R)} \abs{\sigma} = \lim_{k \To \infty} \norm{\sigma}_{L^{\beta_k} (A(R,2R))} 
\leq C \left(\rho(R)^\frac{2}{n-2} R^{-2} \right)^{\frac{n}{2m}} \norm{\sigma}_{L^{m} (A(R/2,5R/2))}.
$$
\end{preuve}

Now, we carry on our study of general Schrödinger operators. We wish to point out a Gagliardo-Nirenberg 
type inequality, which will prove useful later.

\vspace{\baselineskip}


\subsection{The inversion of Schrödinger operators.}

Our purpose is to solve $(\overline{\Delta} + V) \sigma = \tau$, with a convenient $\tau$, and to obtain bounded solutions. 

First, the weighted Sobolev inequality easily yields the 

\begin{lem}\label{estim1}
For $s \geq \frac{2n}{n+2}$, there exists a constant $C(n,s)$ such that
$$
\Forall{\sigma}{C_c^\infty(E)} 
\norm{\sigma}_{L^\frac{ns}{n-2s} (E,\mu_\rho)} \leq  C(n,s) 
\norm{S \overline{\Delta} \sigma}_{L^s ( E ,\rho(r)^\frac{s-1}{n/2-1} vol ) }.
$$ 
\end{lem}

\begin{preuve}
Set $k=\frac{s}{n-2s} \frac{n-2}{2} \geq 1$ and fix $\sigma \in C_c^\infty(E)$.
The weighted Sobolev inequality gives
$$
\frac{1}{S} \norm{\sigma}_{L^\frac{ns}{n-2s} (E,\mu_\rho)}^{2k} \leq \int_M \abs{d (\abs{\sigma}^k)}^{2} = \int_M \abs{\sigma}^k \Delta(\abs{\sigma}^k) 
\leq k \int_M \abs{\sigma}^{2k-1} \Delta \abs{\sigma}, 
$$ 
the Kato inequality then implies
$$
\frac{1}{S} \norm{\sigma}_{L^\frac{ns}{n-2s} (E,\mu_\rho)}^{2k}\leq k \int_{M} \abs{\sigma}^{2k-1} \abs{\overline{\Delta} \sigma}
= k \int_{M} \abs{\sigma}^\frac{n(s-1)}{n-2s} \abs{\overline{\Delta} \sigma}
$$
and the Hölder inequality yields
\begin{eqnarray*}
&\frac{1}{S}& \norm{\sigma}_{L^\frac{ns}{n-2s} (E,\mu_\rho)}^{2k}\\
&\leq& k \left(  \int_M \abs{\overline{\Delta} \sigma}^s \rho(r)^\frac{s-1}{n/2-1} dvol \right)^{1/s}
\left( \int_M \abs{\sigma}^\frac{ns}{n-2s} d\mu_\rho \right)^{1-1/s} 
\end{eqnarray*}
so that eventually
$$
\frac{1}{S} \norm{\sigma}_{L^\frac{ns}{n-2s} (E,\mu_\rho)}
\leq k \left( \int_M \abs{\overline{\Delta} \sigma}^s \rho(r)^\frac{s-1}{n/2-1} dvol \right)^{1/s},
$$
which is indeed the claim.
\end{preuve}

\vspace{\baselineskip}

Now a Nash-Moser iteration yields a $L^\infty$ estimate. 

\begin{lem}\label{estim2}
For every $x>n/2$ and $t\geq 1$, there exists a constant $C(n,x,t)$ 
such that
$$
\Forall{\sigma}{C_c^\infty(E)} 
\norm{\sigma}_{L^\infty (E)}^{\frac{xn}{2x-n} + t} 
\leq C(n,x,t) \norm{S \overline{\Delta} \sigma}_{L^x(E,\rho(r)^\frac{x-1}{n/2-1} vol)}^{\frac{xn}{2x-n}} 
\norm{\sigma}_{L^t(E,\mu_\rho)}^t.
$$ 
\end{lem}

\begin{preuve}
As above, for every $\sigma$ in $C_c^\infty(E)$ and every $k\geq 1$ :
$$
\left( \int_M \abs{\sigma}^{kq} d\mu_\rho \right)^{2/q} \leq k S \int_{M} \abs{\sigma}^{2k-1} \abs{\overline{\Delta} \sigma} dvol. 
$$
Using the Hölder inequality, we deduce:
$$
\left( \int_M \abs{\sigma}^{kq} d\mu_\rho \right)^{2/q} \leq  k S \left(\int_M \abs{\overline{\Delta} \sigma}^x \rho(r)^\frac{x-1}{n/2-1} dvol\right)^{1/x} 
\left( \int_M \abs{\sigma}^{\frac{(2k-1)x}{x-1}}  d\mu_\rho \right)^{1-1/x}. 
$$

Define the sequence $(\beta_i)$ such that $\beta_0 =t$ and
$$
\beta_{i+1} = \frac{q}{2} \left( \frac{x-1}{x} \beta_i +1 \right).
$$
We obtain for every $i \in \Natural$ :
$$
\norm{\sigma}_{L^{\beta_{i+1}}(E,\mu_\rho)}^{\beta_{i+1}} \leq  
\left(q^{-1} \beta_{i+1} S N_x \right)^{\frac{q}{2}}
\left( \norm{\sigma}_{L^{\beta_{i}}(E,\mu_\rho)}^{\beta_{i}} \right)^{\zeta},
$$
where
$$
N_x = \left(\int_M \abs{\overline{\Delta} \sigma}^x \rho(r)^\frac{x-1}{n/2-1} dvol\right)^{1/x} 
$$
and
$$
\zeta = \frac{q(x-1)}{2x} > 1.
$$

Iterating this, we see that for every $i \in \Natural$,
$$
\norm{\sigma}_{L^{\beta_{i}}(E,\mu_\rho)}^{\beta_{i}} \leq 
\left(q^{-1} S N_x  \right)^{\frac{q}{2}\sum_{j=0}^{i-1} \zeta^j}
\left( \prod_{j=1}^{i} \beta_j^{\zeta^{i-j} } \right)^{q/2}
\left( \norm{\sigma}_{L^{\beta_{0}}(E,\mu_\rho)}^{\beta_{0}} \right)^{\zeta^i}.
$$
Thus
$$
\norm{\sigma}_{L^{\beta_{i}}(E,\mu_\rho)} \leq 
\left(q^{-1} S N_x  \right)^{\frac{q}{2 \beta_i} \frac{\zeta^i - 1}{\zeta - 1}}
\left( \prod_{j=1}^{i} \beta_j^{ \zeta^{-j}} \right)^{ \frac{q \zeta^i}{2 \beta_i}}
\left( \norm{\sigma}_{L^{\beta_{0}}(E,\mu_\rho)}^{\beta_{0}} \right)^{\frac{\zeta^i}{\beta_i}}.
$$

Using 
$$
\beta_i = \zeta^i \left( \beta_0 + \frac{q}{2(\zeta-1)} \right) - \frac{q}{2(\zeta-1)},
$$
we see that
$$
\left( \frac{\zeta^i}{\beta_i} \right) \To \frac{1}{\beta_0 + \frac{q}{2(\zeta-1)}}.
$$
Writing 
$$
\log \left( \prod_{j=1}^{i} \beta_j^{ \zeta^{-j}} \right)
= \sum _{j=1}^{i}  j \zeta^{-j} \log \zeta + \sum _{j=1}^{i} \zeta^{-j} \log \frac{\beta_j}{\zeta^j},
$$
we see that this expression has a limit when $i \To \infty$. So
$$
\norm{\sigma}_{L^{\infty}(E,\mu_\rho)} \leq 
\left(q^{-1} S N_x  \right)^{\frac{\frac{q}{2(\zeta-1)}}{\beta_0 + \frac{q}{2(\zeta-1)}}}
\left( \prod_{j=1}^{\infty} \beta_j^{ \zeta^{-j}} \right)^{\frac{\frac{q}{2}}{\beta_0 + \frac{q}{2(\zeta-1)}}}
 \norm{\sigma}_{L^{\beta_{0}}(E,\mu_\rho)}^{\frac{\beta_{0}}{\beta_0 + \frac{q}{2(\zeta-1)}}}.
$$
As
$
\frac{q}{2(\zeta-1)} = \frac{xn}{2x-n}
$
and $\beta_0 =t$, this is what we claimed.
\end{preuve}

\vspace{\baselineskip}

These facts lead to the 

\begin{thm}[Inversion of the Bochner laplacian]\label{inv2}
Choose an element $s$ in $[\frac{2n}{n+2} , \frac{n}{2}[$ and  a number $x > \frac{n}{2}$. Let $\Omega$ be an open set with 
smooth boundary. Then we can define a continuous operator
$$
\overline{\Delta}^{-1} \, : \, L^s (E_\Omega, \rho(r)^\frac{s-1}{n/2-1} vol) \cap L^x (E_\Omega, \rho(r)^\frac{x-1}{n/2-1} vol) 
\To L^\infty (E_\Omega)
$$  
which is an inverse for the Bochner laplacian over $\Omega$, with Dirichlet boundary condition. More precisely, for $\sigma \in C_c^\infty(E_\Omega)$, 
we have the estimate
$$
\norm{\sigma}_{L^\infty (E_\Omega)}^{\frac{s}{n-2s} + \frac{x}{2x-n}} 
\leq C(n,s,x)  \norm{S \overline{\Delta} \sigma}_{L^s(E_\Omega,\rho(r)^\frac{s-1}{n/2-1} vol)}^{\frac{s}{n-2s}} 
\norm{S \overline{\Delta} \sigma}_{L^x(E_\Omega,\rho(r)^\frac{x-1}{n/2-1} vol)}^{\frac{x}{2x-n}}.
$$ 
\end{thm}

\begin{preuve}
The estimate is simply obtained by combining (\ref{estim1}) and (\ref{estim2}).
Given $\psi$ in $C^\infty_c(E_\Omega)$, the classical $L^2$ theory yields a smooth solution 
$\sigma_R$ to the equation $\overline{\Delta} \sigma_R = \psi$ on $\Omega \cap B(R)$, with Dirichlet boundary condition. We extend it into
a continuous function on $\Omega$ by deciding it is zero outside $B(R)$. The $L^\infty$-estimate (which is easily seen to hold for $\sigma_R$,
by looking at the proofs above) gives
$$
\norm{\sigma_{R}}_{L^\infty(E_\Omega)} \leq C(n,s,x)  \norm{S \psi}_{L^s(E_\Omega,\rho(r)^\frac{s-1}{n/2-1} vol)}^{\frac{s}{n-2s}} 
\norm{S \psi}_{L^x(E_\Omega,\rho(r)^\frac{x-1}{n/2-1} vol)}^{\frac{x}{2x-n}}.
$$
For every compact set $K$, there is an $R_K$ such that the family 
$(\sigma_R\restric{K}, \, R\geq R_K)$ is uniformly bounded in $C^\infty(E_K)$ (by elliptic regularity), so that Ascoli yields a sequence 
converging in $C^\infty(E_K)$. By a diagonal extraction, we find a sequence $(\sigma_{R_i})$ which converges to 
$\sigma$ in $C^0_c(E_\Omega)$. $\sigma$ is easily seen to be a weak solution of $\overline{\Delta} \sigma = \psi$, it is therefore smooth and thus a strong solution. For every compact set $K$, 
\begin{eqnarray*}
\norm{\sigma}_{L^\infty(E_K)} &=& \lim_{i \To \infty} \norm{\sigma_{R_i}}_{L^\infty(E_K)} \\
&\leq& C(n,s,x)  \norm{S \psi}_{L^s(E_\Omega,\rho(r)^\frac{s-1}{n/2-1} vol)}^{\frac{s}{n-2s}} 
\norm{S \psi}_{L^x(E_\Omega,\rho(r)^\frac{x-1}{n/2-1} vol)}^{\frac{x}{2x-n}},
\end{eqnarray*}
hence a $L^\infty$-estimate on $\Omega$. We can thus define an operator $\overline{\Delta}^{-1}$ on $C^\infty_c(E_\Omega)$ which is continuous for the expected norms. We then extend it by continuity.
\end{preuve}

\vspace{\baselineskip}

By a perturbation technique, we deduce an analogous result for Schrödinger operators. 

\begin{thm}[Inversion of Schrödinger operators]\label{inv3}
Set $\frac{2n}{n+2} \leq s < \frac{n}{2}$ and $x > \frac{n}{2}$. Then there exists a positive number $\eta(n,s,x,S)$ such that, 
given an open set with smooth boundary $\Omega$ and a potential $V$ satisfying
$$
S \max(\norm{V_-}_{L^s(\Omega,\rho(r)^\frac{s-1}{n/2-1} vol)},\norm{V_-}_{L^x(\Omega,\rho(r)^\frac{x-1}{n/2-1} vol)}) < \eta(n,s,x),
$$
there is a continuous operator
$$
(\overline{\Delta} + V)^{-1} \, : \, L^s (E_\Omega, \rho(r)^\frac{s-1}{n/2-1} vol) \cap L^x (E_\Omega, \rho(r)^\frac{x-1}{n/2-1} vol) 
\To L^\infty (E_\Omega). 
$$
\end{thm}

\begin{preuve}
First, the previous analysis works for $H:=\overline{\Delta} + V_+$ as well as for $\overline{\Delta}$.
Then define $\eta(n,s,x)$ to be $S$ divided by the norm of
$$
H^{-1} \, : \, L^s (E_\Omega, \rho(r)^\frac{s-1}{n/2-1} vol) \cap L^x (E_\Omega, \rho(r)^\frac{x-1}{n/2-1} vol) \To L^\infty (E_\Omega),
$$
so that, under our assumption, 
$$
V_-\, : \, L^\infty (E_\Omega)
\To L^s (E_\Omega, \rho(r)^\frac{s-1}{n/2-1} vol) \cap L^x (E_\Omega, \rho(r)^\frac{x-1}{n/2-1} vol),
$$ 
is a continuous operator whose norm is strictly inferior to $\eta(n,s,x)/S$ and $H^{-1} V_-$  is a continuous endomorphism
of $L^\infty (E_\Omega)$, with norm strictly inferior to $1$. The operator $\Id + H^{-1} V_-$ 
is then an automorphism of $L^\infty (E_\Omega)$. So we can define the continuous operator 
$(\Id + H^{-1} V_-)^{-1} H^{-1} = (\overline{\Delta} + V)^{-1}$, from  
$L^s (E_\Omega, \rho(r)^\frac{s-1}{n/2-1} vol) \cap L^x (E_\Omega, \rho(r)^\frac{x-1}{n/2-1} vol)$ to $L^\infty (E_\Omega)$. 
\end{preuve}


\section{Applications.}

\subsection{$L^2$-cohomology.}

Our study of Schrödinger operators gives geometric information as soon as the potential depends only on the curvature tensor. 
For instance, if the weighted Sobolev inequality is true, the vanishing theorem (\ref{annulation}) forces the kernel of such 
"geometric operators" to be trivial, under integral assumptions on the curvature. We discuss here the case of the Hodge laplacian 
$\Delta = d d^* + d^* d$. It is well known that this operator, when acting on $k$-forms, admits the Weitzenböck decomposition
$$
\Delta^k = \overline{\Delta} + \mathcal{R}^k,
$$
where $\mathcal{R}^k$ is a field of symmetric endomorphisms of the vector bundle of $k$-exterior forms, depending only 
on the curvature. In particular, $\mathcal{R}^1 = \Ric$. Our results apply and we can obtain information on the (reduced) 
$L^2$-cohomology $\mathcal{H}_{L^2}(M)$. We refer to \cite{Car1} for the definitions. The point is that $\mathcal{H}_{L^2}^k(M)$ can 
be identified with the kernel of $\Delta^k$, seen as an unbounded operator on $L^2$ $k$-forms. 

We can indeed generalize G. Carron's results in \cite{Car1}. Before stating our theorem, we need to introduce the following decreasing 
function, derived from the Euler $\Gamma$ function ($q=2n/(n-2)$):
$$
\iota_q \, : \, x \, \mapsto 
\frac{2}{x}  \left( \frac{\Gamma \left(\frac{x+q}{2} \right)}{\Gamma\left(\frac{x}{2}\right)} \right)^{2/q}.
$$

\begin{thm}[$L^2$-cohomology]
Let $M^n$, $n \geq 3$, be a connected complete riemannian manifold with nonnegative Ricci curvature. Assume that there
exists $o \in M$, $\nu > 2$ and  $C_o >0$ such that 
$$
\forall R_2 \geq R_1 > 0, \  \frac{V(o,R_2)}{V(o,R_1)} \geq C_o \left( \frac{R_2}{R_1} \right)^\nu.
$$
Then $\mathcal{H}_{L^2}^1(M) = \set{0}$. Let $k \geq 2$ be an integer. 
\begin{itemize}
\item If $\norm{\mathcal{R}_-^k}_{L^{n/2}(\rho_o(r_o) vol)} < \infty$,
then $\dim \mathcal{H}_{L^2}^k(M) < \infty$. \\
\item If $S_{o}(M) \norm{\mathcal{R}_-^k}_{L^{n/2}(\rho_o(r_o) vol)} < 1$,
then $\mathcal{H}_{L^2}^k(M) = \set{0}$. \\
\item If, for some integer $N_0 \geq \binom{n}{k}$, $S_{o}(M) \norm{\mathcal{R}_-^k}_{L^{n/2}(\rho_o(r_o) vol)} \leq  \frac{\iota_q(k)}{\iota_q(N_0)}$, then 
$\dim \mathcal{H}_{L^2}^k(M) \leq N_0.$\\
\item Given $\frac{2n}{n+2} \leq s < n/2 < x$, there exists a constant $C = C(n,s,x)$ such that 
the dimension of $\mathcal{H}_{L^2}^k(M)$ is bounded by 
$$ 
\binom{n}{k} \, \max \left( 1 , C \norm{S_{o}(M) \mathcal{R}_-^k}_{L^x(E,\rho_o(r_o)^\frac{x-1}{n/2-1} vol)}^{x} 
\norm{S_{o}(M) \mathcal{R}_-^k}_{L^s(E,\rho_o(r_o)^\frac{s-1}{n/2-1} vol)}^{\frac{s(2x-n)}{n-2s}} \right).
$$
\end{itemize}
\end{thm}

\begin{cor}
Let $M^n$, $n \geq 3$, be a connected complete riemannian manifold with nonnegative Ricci curvature. Assume that there
exists $o \in M$, $\nu > 2$ and  $C_o >0$ such that 
$$
\forall R_2 \geq R_1 > 0, \  \frac{V(o,R_2)}{V(o,R_1)} \geq C_o \left( \frac{R_2}{R_1} \right)^\nu
$$
and the Riemann curvature tensor satisfies
$$
\left(\int_M \abs{R}^{\frac{n}{2}} \rho_o(r_o) dvol \right)^{\frac{2}{n}} < \infty.
$$
Then the $L^2$-cohomology of $M$ is finite dimensional. 
\end{cor}

We omit the proof, which consists in using the weighted Sobolev inequality (\ref{inegalite}), in order to make the techniques of \cite{Car1}
work. The vanishing results stem from \ref{annulation}, of course.

\vspace{\baselineskip}


\subsection{Ricci flat manifolds.}\label{ricciflat}

\subsubsection{Flatness criterions.}

We want to explain here why the weighted Sobolev and Hardy inequalities help understanding Ricci flat manifolds. In particular, 
they emphasize some rigidity properties of these manifolds, under volume growth assumptions. We will show that if their curvature is small, 
in some integral sense, then they are actually flat. 

The key tool is a property of the Weyl tensor $W$ of a Ricci-flat manifold with dimension $n \geq 4$ : it obeys the nonlinear equation 
$$
\overline{\Delta} W = W * W,
$$
where the right-hand side is a quadratic expression in the Weyl curvature \cite{Bes}. In particular, $W$ is either identically zero, or 
vanishes only on a set of zero measure. So, outside a set of zero measure, $\abs{W}$ is smooth and satisfies the estimate
$$
\abs{\Delta \abs{W}} \leq c(n) \abs{W}^2,
$$
where $c(n)$ is a universal constant, depending only on the dimension $n$.

For every $k \geq 1$, 
$$
\Delta \abs{W}^k = k \abs{W}^{k-1} \Delta \abs{W} - k(k-1) \abs{W}^{k-2} \abs{d\abs{W}}^2 \leq k c(n) \abs{W}^{k+1}.
$$
It turns out that this inequality is still true for some $k<1$. This is made possible by the refined Kato inequality (\cite{BKN}, \cite{CGH}), 
which says that the Weyl tensor $W$ of a Ricci-flat $n$-manifold satisfies almost everywhere
$$
\abs{d\abs{W}}^2 \leq \frac{n-1}{n+1} \abs{\nabla W}^2.
$$

>From this, one can deduce that almost everywhere, 
$$
\Delta \abs{W}^\gamma \leq c(n) \gamma \abs{W}^{1+\gamma},
$$
with 
$$
\gamma := \frac{n-3}{n-1}.
$$
Indeed, note $\frac{n-1}{n+1}= \frac{1}{2-\gamma}$ and write 
\begin{eqnarray*}
\Delta \abs{W}^\gamma &=& \gamma \abs{W}^{\gamma-1} \Delta \abs{W} + \gamma(1-\gamma) \abs{W}^{\gamma-2} \abs{d\abs{W}}^2\\
&=& \gamma \abs{W}^{\gamma-2} \left( \frac{1}{2} \Delta\abs{W}^2 + \abs{d\abs{W}}^2\right) + 
\gamma(1-\gamma) \abs{W}^{\gamma-2} \abs{d\abs{W}}^2\\
&=& \gamma \abs{W}^{\gamma-2} \left((W,\overline{\Delta}W) - \abs{\nabla W}^2 \right) + \gamma(2-\gamma) \abs{W}^{\gamma-2} \abs{d\abs{W}}^2\\
&\leq& c(n)\gamma \abs{W}^{\gamma+1} - \gamma\abs{W}^{\gamma-2} \abs{\nabla W}^2 + \gamma\abs{W}^{\gamma-2} \abs{\nabla W}^2\\
&=& c(n)\gamma \abs{W}^{\gamma+1}. 
\end{eqnarray*}

Now, given $k\geq \gamma$, we can write $k=\gamma l$, $l\geq 1$ and then
\begin{eqnarray*}
\Delta \abs{W}^k &=& \Delta (\abs{W}^\gamma)^l\\
&=& l (\abs{W}^\gamma)^{l-1} \Delta (\abs{W}^\gamma) - l(l-1) (\abs{W}^\gamma)^{l-2} \abs{d(\abs{W}^\gamma)}^2 \\
&\leq& l (\abs{W}^\gamma)^{l-1} c(n) \gamma  \abs{W}^{\gamma+1}\\
&=& k c(n) \abs{W}^{k+1}.
\end{eqnarray*}

Thus 
\begin{equation}\label{weyl}
\Delta \abs{W}^k \leq c(n) k \abs{W}^{1+k},
\end{equation}
is true for any $k\geq \gamma$. 

\vspace{\baselineskip}

With this differential inequality in hand, we can prove flatness and curvature decay results. To express them, we need the
\begin{defn}
The "Sobolev-curvature" invariant of a connected complete manifold $M^n$ is defined by 
$$
\mathcal{SC} (M):= \inf_{o \in M} \left[ S_{o}(M) \left(\int_M \abs{R}^{\frac{n}{2}}
\rho_o(r_o) dvol \right)^{\frac{2}{n}} \right],
$$
where $R$ is the Riemann curvature tensor. We also define a "Hardy-curvature" invariant :
$$
\mathcal{HC} (M):= \inf_{o \in M} \left[ H_{o}(M)^2 \sup_M (\abs{R} r_o^2) \right],
$$
We use the convention $0.\infty=\infty$.
\end{defn}

Now, let us phrase one of our main results.
\begin{thm}[Flatness criterion (1)]\label{plate}
We consider a connected complete Ricci-flat manifold $M^n$, with $n \geq 4$. Assume $\mathcal{SC} (M) < \frac{4}{n c(n)}$.
Then $M$ is flat. 
\end{thm}

\begin{preuve}
It is a consequence of the vanishing theorem \ref{annulation}, applied to the operator $\Delta - c(n) \abs{W}$ and 
the section $\abs{W}$, thanks to our weighted Sobolev inequality. Setting $m=\frac{n}{2}$ in (\ref{annulation}), we obtain $W=0$, 
and as $\Ric=0$, $M$ is flat.  
\end{preuve}

\begin{rem}
The "threshold" value $\frac{4}{nc(n)}$ can definitely be improved : in general, it can be replaced by $\frac{4}{\max(n-2,4\gamma)c(n)}$ ; and if $M$ satisfies (\ref{volboules2}), it can even be replaced to $\frac{4}{\max(\nu-2,4\gamma)c(n)}$. The idea is to use the Hölder inequality and an upper bound on the volume growth, so as to estimate the suitable integral. 
\end{rem}

\begin{cor}
Let $M^n$, $n \geq 4$, be a connected complete Ricci-flat manifold. Assume there
exists $o$ in $M$, $\nu > 2$ and  $C_o >0$ such that 
$$
\forall R_2 \geq R_1 > 0, \  \frac{V(o,R_2)}{V(o,R_1)} \geq C_o \left( \frac{R_2}{R_1} \right)^\nu.
$$
Then there is a constant $\epsilon(n,C_o,\nu)$ such that $M$ is flat as soon as 
$$
\int_M \abs{W}^{\frac{n}{2}} \rho_o(r_o) dvol < \epsilon(n,C_o,\nu).
$$
\end{cor}

\vspace{\baselineskip}

There is also a flatness criterion based on the Hardy inequality (\ref{hardy}).

\begin{thm}[Flatness criterion (2)]\label{plate2}
We consider a connected complete Ricci-flat manifold $M^n$, with $n \geq 4$. There is a constant $\epsilon(n)$ such that if 
$$
\mathcal{HC} (M) < \epsilon(n),
$$
then $M$ is flat. 
\end{thm}  

\begin{demo}
Choose $o$ in $M$ such that $ H_{o}(M)^2 \sup_M (\abs{W} r_o^2) < \epsilon(n)$ ($\epsilon(n)$ will be determined at the end of the proof)
and set $H=H_{o}(M)$, $K=\sup_M (\abs{W} r_o^2)$. We consider, for large $R$, a smooth function $\chi$ which is equal to $1$ on $B(R)$, equal 
to $0$ on $M \backslash B(2R)$, has values in $[0,1]$ and satisfies $\displaystyle{\abs{\nabla \chi} \leq \frac{2}{R}}$. We also work 
with a number $k \geq 5/4$, which will be fixed later. The Hardy inequality says that
$$
\int_M \chi^2 \abs{W}^{2k} r^{-1} \leq H \int_M \abs{d(\chi^2 \abs{W}^{2k})}.
$$
The right hand side can be bounded via triangle and Cauchy-Schwarz inequalities :
$$
\int_M \abs{d(\chi^2 \abs{W}^{2k})} \leq 2 \int_M \chi \abs{d\chi} \abs{W}^{2k} 
+ 2 \left( \int_M \chi^2 \abs{W}^{2k} r^{-1}\right)^{1/2} \left( \int_M \chi^2 \abs{d(\abs{W}^{k})}^2 r \right)^{1/2}.
$$
Set $k':= k -1/4$. So as to perform integration by parts, we kill the $r$ in the lattest integral :
$$
\int_M \chi^2 \abs{d(\abs{W}^{k})}^2 r = (k/k')^2 \int_M \chi^2 \abs{d(\abs{W}^{k'})}^2 \abs{W}^{1/2} r 
\leq (k/k')^2 K^{1/2} \int_M \chi^2 \abs{d(\abs{W}^{k'})}^2.
$$
Integration by parts and (\ref{weyl}) yield : 
\begin{eqnarray*}
& &\left( \int_M \abs{\chi d(\abs{W}^{k'})}^2 \right)^{1/2} = \left( \int_M \abs{d(\chi \abs{W}^{k'}) - \abs{W}^{k'} d\chi }^2\right)^{1/2} \\
& &\leq   \left(\int_M \abs{W}^{2 k'} \abs{d\chi }^2 \right)^{1/2} + \left( \int_M \abs{d(\chi \abs{W}^{k'})}^2 \right)^{1/2} \\
& &\leq  2 \left( \int_M \abs{W}^{2 k'} \abs{d\chi }^2 \right)^{1/2} +  \left(\int_M \chi^2 \abs{W}^{k'} \Delta \abs{W}^{k'}\right)^{1/2} \\
& &\leq  2 \left( \int_M \abs{W}^{2 k'} \abs{d\chi }^2 \right)^{1/2} +  k'^{1/2} c(n)^{1/2} \left(\int_M \chi^2 \abs{W}^{2k'+1} \right)^{1/2}\\
& &\leq  2 \left(\int_M \abs{W}^{2 k -1/2} \abs{d\chi }^2 \right)^{1/2} 
+ ((k-1/4) c(n))^{1/2} K^{1/4} \left(\int_M \chi^2 \abs{W}^{2k} r^{-1} \right)^{1/2}. 
\end{eqnarray*}
Thus
\begin{eqnarray*}
\int_M \chi^2 \abs{W}^{2k} r^{-1} &\leq& 2 H \int_M \chi \abs{d\chi} \abs{W}^{2k} \\
&+& 16 H K^{1/4} k/(4k-1)   \left( \int_M \chi^2 \abs{W}^{2k} r^{-1} \right)^{1/2} \left(\int_M \abs{W}^{2 k -1/2} \abs{d\chi }^2 \right)^{1/2} \\
&+&  4H (c(n) K)^{1/2} k/(4k-1)^{1/2}  \int_M \chi^2 \abs{W}^{2k} r^{-1}. 
\end{eqnarray*}
We want to ensure $4H (c(n) K)^{1/2} k/(4k-1)^{1/2} $ is strictly less than $1$. Indeed, this is
equivalent to $16 H^2 K c(n) k^2 -4k +1 < 0$ ; this trinomial has two positive roots for $H^2 K$ (and thus $\epsilon(n)$) small enough, 
and we can choose $k$ to be half the sum of theses roots : $k := (8 H^2 K c(n))^{-1}$. Then we obtain
\begin{eqnarray*}
& &(1 - 4H (c(n) K)^{1/2} k/(4k-1)^{1/2} ) \int_{B(R)} \abs{W}^{2k} r^{-1} \\
& &\leq  4 H R^{-1} \int_M \abs{W}^{2k} 
 +  32 H  K^{1/4} k/(4k-1) R^{-1} \left( \int_M \abs{W}^{2k} r^{-1} \right)^{1/2} \left(\int_M \abs{W}^{2 k -1/2} \right)^{1/2}.
\end{eqnarray*}
If we choose $\epsilon(n)$ small enough (so that $k$ is big enough), the integrals on the right hand side are finite (recall $W$ has quadratic
decay and the volume growth is at most euclidian) ; one can eventually take $\epsilon(n) = \frac{1}{2 (n+1) c(n)}$. Letting 
$R \To \infty$, we find $W=0$, and with $Ric = 0$, $M$ is flat. 
\end{demo}

\begin{cor}
Let $M^n$, $n \geq 4$, be a connected complete Ricci-flat manifold. Assume there
exists $o$ in $M$, $\nu > 1$ and  $C_o >0$ such that 
$$
\forall R_2 \geq R_1 > 0, \  \frac{V(o,R_2)}{V(o,R_1)} \geq C_o \left( \frac{R_2}{R_1} \right)^\nu.
$$
Then there is a constant $\epsilon(n,C_o,\nu)$ such that $M$ is flat as soon as 
$$
\sup_M (\abs{W} r_o^2) < \epsilon(n,C_o,\nu).
$$
\end{cor}


\subsubsection{Curvature decay.}

In the preceding paragraph, we have seen that when $\mathcal{SC}(M)$ is small, the curvature vanishes. 
Now using the decay lemmas of \ref{decaylemmas}, we can show that if $\mathcal{SC}(M)$ is only finite, then the curvature decays at infinity.   

We first prove the
\begin{prop}\label{deccourb1}
We consider a connected complete Ricci-flat manifold $M^n$, $n \geq 4$ such that $\mathcal{SC}(M)$ is finite.
Then for any point $o$ in $M$,
$$
\sup_{S(o,R)} \abs{W} = o(R^{- 2}).
$$
\end{prop}

Before proving this, let us state a consequence of our work.
\begin{cor}
We consider a connected complete Ricci-flat manifold $M^n$ with $n \geq 4$. Assume there exists $o$ in $M$, 
$\nu > 2$ and  $C_o >0$ such that 
$$
\forall R_2 \geq R_1 > 0, \  \frac{V(o,R_2)}{V(o,R_1)} \geq C_o \left( \frac{R_2}{R_1} \right)^\nu
$$
and the curvature satisfies 
$$
\int_M \abs{W}^{\frac{n}{2}} \rho_o(r_o) dvol < + \infty.
$$
Then
$$
\sup_{S(o,R)} \abs{W} = o(R^{- 2}).
$$
\end{cor}

\begin{rem}
This should be compared with the result of \cite{CT} : $\sup_{S(o,R)} \abs{W} = O(R^{- 2})$ as soon as $n=4$, $\Ric=0$ and $W \in L^2$. 
\end{rem}

\begin{rem}
If we assume $W$ behaves like $r^{-\sigma}$, the assumption 
$$
\int_M \abs{W}^{\frac{n}{2}} \rho_o(r_o) dvol < + \infty
$$ 
is equivalent to
$\sigma>2$ : the above result therefore turns an integral estimate into the pointwise estimate we can hope. The next theorem will point out 
an automatic improvement of the decay ; it is another rigidity phenomenon.
\end{rem}

\vspace{\baselineskip}

\begin{demo}
As $(\Delta - c(n) \abs{W}) \abs{W} \leq 0$, we want to apply the lemma \ref{lemdec3} with the operator $\Delta - c(n) \abs{W}$. To do this, 
we set $x = nq/4$, so that in particular $x-n/2 = \frac{n}{n-2}$ and use lemma \ref{lemdec1} with $m=n/2$, which implies that for large $R$ : 
\begin{eqnarray*}
&&\left( \int_{A(R,2R)} \abs{W}^x \rho(r)^\frac{x-1}{n/2-1} dvol \right)^\frac{1}{x-n/2}  \\
&\leq& C \left( \rho(R)^\frac{2x}{n-2} \int_{A(R,2R)} \abs{W}^{nq/4} d\mu_\rho \right)^\frac{1}{x-n/2}\\
&=& C \rho(R)^\frac{n}{n-2} \left( \int_{A(R,2R)} \abs{W}^{nq/4} d\mu_\rho \right)^\frac{n-2}{n}\\
&\leq& C \rho(R)^\frac{n}{n-2} R^{-2}  \int_{A(R/2,R)} \abs{W}^{n/2} dvol \\
&\leq& C \rho(R)^\frac{2}{n-2} R^{-2}  \int_{A(R/2,R)} \abs{W}^{n/2} \rho(r) dvol \\
&\leq& C \rho(R)^\frac{2}{n-2} R^{-2}.
\end{eqnarray*}  

We can use lemma \ref{lemdec3} with $m=n/2$ : 
\begin{eqnarray*}
\sup_{S(R)} \abs{W}^{n/2} 
&\leq& C \left( \rho(R)^\frac{2}{n-2} R^{-2} \right)^{n/2} \int_{A(R/2,5R/2)} \abs{W}^{n/2} d\mu_\rho \\
&\leq& C R^{-n} \int_{A(R/2,5R/2)} \abs{W}^{n/2} \rho(r) dvol.
\end{eqnarray*}
Hence $W = o(r^{-2})$ (since the right-hand side integral tends to zero).
\end{demo}

\vspace{\baselineskip}

In general, such a quadratic curvature decay is not so meaningful : actually, any smooth connected noncompact manifold admits a metric 
with quadratic curvature decay \cite{Gro}. Note however that a riemannian manifold with nonnegative Ricci curvature, maximal volume growth 
($\nu=n$) and quadratic curvature decay has finite topological type \cite{SS}. In case the volume growth is not maximal, such a strong
topological consequence is not known. 

We would like to point out a consequence of the quadratic curvature decay. Applying one of the results of \cite{LS}, it yields the
\begin{cor}
Let $M^n$, $n \geq 4$, be a connected complete Ricci-flat manifold. Assume there exists $o$ in $M$, 
$\nu > 2$ and  $C_o >0$ such that 
$$
\forall R_2 \geq R_1 > 0, \  \frac{V(o,R_2)}{V(o,R_1)} \geq C_o \left( \frac{R_2}{R_1} \right)^\nu,
$$
$$
V(o,R) = o(R^n)
$$
and the curvature satisfies 
$$
\int_M \abs{W}^{\frac{n}{2}} \rho_o(r_o) dvol < + \infty.
$$
Then the integral of the Chern-Gauss-Bonnet form is an integer. 
\end{cor} 

\begin{rem}
If $n=4$, this means 
$$
\frac{1}{8\pi^2} \int_M \abs{W}^2 dvol \in \Integer.
$$
In particular, if $\int_M \abs{W}^2 dvol < 8\pi^2$, $M$ is flat.
\end{rem}

Now, it is well known that manifolds with faster than quadratic curvature decay enjoy nice properties \cite{Abr}. This motivates 
our quest for a better estimate on the curvature. The key is the refined Kato inequality.

\begin{thm}[Curvature decay (1)]\label{deccourb2}
We consider a connected complete Ricci-flat manifold $M^n$, $n \geq 4$, such that $\mathcal{SC}(M)$ is finite.  
Fix a point $o$ in $M$ and assume there exists $\nu > 2$ and $A_o>0$ such that 
$$
\forall R\geq 1, \,  V(o,R) \geq A_o R^\nu.
$$
Then 
$$
\sup_{S(o,R)} \abs{W} = O(R^{- b})
$$ 
for $b=2$ and every $b < \frac{\nu-2}{\gamma} = \frac{(\nu-2)(n-1)}{n-3}$.
\end{thm}

\begin{demo}
Set $w=\abs{W}^\gamma$ and 
$
b_0 = \sup \set{b>0 \, / \, w = O\left( \left[r^2/V(r)\right]^{b} \right)}.
$
We know, from the previous proposition, that $w = O(r^{-2\gamma})$ ; since $V(r) \leq \omega_n r^n$ (Bishop), this implies
$w = O(V(r)^{-2\gamma/n}) = O([r^2 / V(r)]^{2\gamma/n})$, so that $b_0$ is a positive number. 
Suppose $b_0 < 1$. We can choose $b_1 > 0$, $m>0$ such that  
$$
m> \frac{n}{b_1(n-2)} > \frac{n}{b_0(n-2)} > \frac{n}{n-2}.
$$ 
Since $b_1 < b_0$, $w = O\left( \left[r^2/V(r)\right]^{b_1} \right)$, so that for any $R>0$, 
\begin{eqnarray*}
\int_{A(R,2R)} \abs{w}^{m} d\mu_\rho &\leq& C \left[R^2 /V(R)\right]^{m b_1} \rho(R)^{-\frac{2}{n-2}} V(R)\\ 
&=& C \left[R^2 /V(R)\right]^{m b_1 - \frac{n}{n-2}} \\
&\leq& C R^{-(\nu-2)\left(m b_1 - \frac{n}{n-2}\right)}.
\end{eqnarray*}
This implies
$$
\int_M \abs{w}^{m} d\mu_\rho < + \infty.
$$
Recall that almost everywhere
$$
(\Delta - \gamma c(n)\abs{W}) w \leq 0.
$$
We intend to apply lemma \ref{lemdec2} to the function $w$, which is unfortunately not locally Lipschitz. To overcome this, once again, 
we consider $u_\epsilon := \sqrt{\abs{W}^2 +\epsilon}$, $\epsilon>0$. Direct computation yields almost everywhere
\begin{eqnarray*}
u_\epsilon^\gamma \Delta u_\epsilon^\gamma &=& \gamma u_\epsilon^{2 \gamma -2} \left( \abs{W} \Delta \abs{W} 
- \epsilon u_\epsilon^{-2} \abs{d\abs{W}}^2 \right) + \gamma (1-\gamma) u_\epsilon^{2\gamma-4} \abs{W}^2 \abs{d\abs{W}}^2 \\ 
&\leq& \gamma u_\epsilon^{2 \gamma -2} \left( \abs{W} \Delta \abs{W} +
(1-\gamma) \abs{d\abs{W}}^2 \right)
\end{eqnarray*}
and, using the refined Kato inequality as in the proof of (\ref{weyl}), we find (everywhere)
$$
u_\epsilon^\gamma \Delta u_\epsilon^\gamma \leq \gamma u_\epsilon^{2 \gamma} ( W, \overline{\Delta} W ).
$$   
As in the proof of (\ref{annulation}), by making $\epsilon$ go to zero, we are able to obtain the first inequality in the proof of 
lemma \ref{lemdec2} ($m>\frac{n}{n-2}$). Eventually,
$$
\int_{M \backslash B(R)} \abs{w}^{m} d\mu_\rho =O(R^{-a}),
$$
for some $a>0$ which is independent of the choice of $m$ in a neighbourhood of $\frac{n}{b_0(n-2)}$. 
Now, applying the lemma \ref{lemdec3} (again, one must adapt the proof because $w$ is not locally Lipschitz), with this $m$, we find for large $R$ :
\begin{eqnarray*}
\sup_{S(R)} w 
&\leq& C \left(\rho(R)^\frac{2}{n-2} R^{-2} \right)^{\frac{n}{2m}} R^{-a/m} \\
&=& C \left[ R^{2} / V(R) \right]^\frac{n}{m(n-2)} R^{-a/m}\\
&\leq& C \left[ R^{2} / V(R) \right]^{\frac{n}{m(n-2)}+\frac{a}{nm}},
\end{eqnarray*}
where we again used the euclidian upper bound on the volume growth of balls.
When $m$ goes to $\frac{n}{b_0(n-2)}$, the exponent tends to $b_o + \frac{b_o (n-2) a}{n^2}$ : if we choose $m$ sufficiently close to
$\frac{n}{b_0(n-2)}$, we obtain a contradiction to the definition of $b_o$. So $b_o  \geq 1$ and, with the lower bound on the volume growth, we are done. 
\end{demo}

\begin{cor}
We consider a connected complete Ricci-flat manifold $M^n$ with $n \geq 4$. Assume there exists $o$ in $M$, 
$\nu > 2$ and  $C_o >0$ such that 
$$
\forall R_2 \geq R_1 > 0, \  \frac{V(o,R_2)}{V(o,R_1)} \geq C_o \left( \frac{R_2}{R_1} \right)^\nu
$$
and the curvature satisfies 
$$
\int_M \abs{W}^{\frac{n}{2}} \rho_o(r_o) dvol < + \infty.
$$
Then 
$$
\sup_{S(o,R)} \abs{W} = O(R^{- b})
$$ 
for $b=2$ and every $b < \frac{\nu-2}{\gamma} = \frac{(\nu-2)(n-1)}{n-3}$.
\end{cor}

\vspace{\baselineskip}

Let us point out the topological consequence we were expecting.

\begin{cor}[Finite topology]
A connected complete Ricci-flat manifold $M^n$, $n \geq 4$, for which there exists a point $o$, 
$\nu > 4\frac{n-2}{n-1}$ and  $C_o >0$ such that 
$$
\forall R_2 \geq R_1 > 0, \  \frac{V(o,R_2)}{V(o,R_1)} \geq C_o \left( \frac{R_2}{R_1} \right)^\nu
$$
and whose curvature satisfies 
$$
\int_M \abs{W}^{\frac{n}{2}} \rho_o(r_o) dvol < + \infty.
$$
is homeomorphic to the interior of a compact manifold with boundary. More precisely, there is a connected open subset $U$ of $M$
which has compact closure, smooth boundary and such that $M\backslash U$ is a connected manifold with boundary
which is diffeomorphic to $N \times \Real_+$ for some closed connected $n-1$-manifold $N$. Furthermore, if $V(o,R) \asymp R^\nu$ 
with $4\frac{n-2}{n-1}<\nu<n$, we know that $N$ either has trivial tangent bundle or infinite fundamental group ; in case $2\frac{3n-7}{n-1} < \nu < n$ 
(faster than quartic curvature decay and strictly subeuclidian volume growth), $N$ always has infinite fundamental group.   
\end{cor}

\begin{preuve}
The theorem implies $M$ has faster than quadratic curvature decay so that \cite{Abr},\cite{GPZ} apply.
\end{preuve}

\vspace{\baselineskip}

One can wonder whether the limiting decay exponent is indeed attained. Actually, this is the case.

\begin{thm}[Curvature decay (2)]\label{deccourb3}
We consider a connected complete Ricci-flat manifold $M^n$, $n \geq 4$, such that $\mathcal{SC}(M)$ is finite.  
Fix a point $o$ in $M$ and assume there exists $\nu > 4\frac{n-2}{n-1}$ and $A_o>0$ such that 
$$
\forall R\geq 1, \,  V(o,R) \geq A_o R^\nu.
$$
Then 
$$
\sup_{S(o,R)} \abs{W} = O\left(r^{- \frac{(\nu-2)(n-1)}{n-3}}\right).
$$ 
\end{thm}

\begin{demo}
In \cite{Gur}, Gursky introduced the following operator :
$$
L_g := \Delta_g +  \frac{n-2}{4(n-1)} Scal_g - \gamma c(n) \abs{W}_g.
$$
It turns out that this operator is conformally invariant in the following sense : if $\phi$ is a smooth positive function, 
\begin{equation}\label{invconf}
L_{\phi^{\frac{4}{n-2}} g} = \phi^{-\frac{n+2}{n-2}} L_g (\phi .).
\end{equation}
We intend to use this property to find in the conformal class of $g$ a new metric $\tilde{g}$ such that outside a compact set 
$$
L_{\tilde{g}} = \Delta_{\tilde{g}},
$$
i.e.
$$
\frac{n-2}{4(n-1)} Scal_{\tilde{g}} - \gamma c(n) \abs{W}_{\tilde{g}} = 0.
$$
We seek $\tilde{g}$ in the form of $\tilde{g} = (1+u)^{\frac{4}{n-2}} g$, where $g$ is our Ricci-flat metric and $u$ is a smooth function 
to determine. Applying (\ref{invconf}) to the constant function $1$, we find
$$
L_{\tilde{g}} (1) = L_{(1+u)^{\frac{4}{n-2}} g}(1) = (1+u)^{-\frac{n+2}{n-2}} L_g (1+u),
$$
so that, since $Scal_g = 0$,
$$
\frac{n-2}{4(n-1)} Scal_{\tilde{g}} - \gamma c(n) \abs{W}_{\tilde{g}} = (1+u)^{-\frac{n+2}{n-2}} ( \Delta_g u - \gamma c(n) \abs{W}_g) (1+u).
$$
We thus have to solve
\begin{equation}\label{edp}
\Delta_g u - \gamma c(n) \abs{W}_g u = \gamma c(n) \abs{W}_g.
\end{equation}
Let us solve it on $M \backslash B_g(o,R)$, with large $R$ (let us assume $S_g(o,R)$ is smooth, this not a problem). We would like to use the inversion theorem 
\ref{inv3} with $\Delta_g - \gamma c(n) \abs{W}_g$. The assumption $\nu > 4\frac{n-2}{n-1}$ ensures $\frac{\nu-2}{\gamma} > 2$ 
: theorem \ref{deccourb2} says that $\abs{W} = O(r^{-b})$ for some $b>2$. In particular, using Bishop's upper bound on the volume growth, one sees that for small $\delta>0$,
$$
\int_M \abs{W}^{n/2 - \delta} \rho(r)^\frac{n/2 - \delta - 1}{n/2-1} dvol < \infty. 
$$ 
Choosing $R$ sufficiently large, we ensure 
$$
S_{o}(M) \int_{M \backslash B_g(o,R)} \abs{W}^{n/2 \pm \delta} \rho(r)^\frac{n/2 \pm \delta - 1}{n/2-1} dvol < \eta(n/2,n/2-\delta,n/2+\delta). 
$$ 
So \ref{inv3} yields a bounded solution $u$ of (\ref{edp}) on $M \backslash B_g(o,R)$, and by enlarging $R$ if necessary, we can assume 
$\norm{u}_{L^\infty} < 1$. Extending $u$ to the whole $M$ in a convenient way, we obtain a metric $\tilde{g}$ which is conformally 
quasi-isometric to $g$ and such that its Gursky operator and its laplacian coincide outside some ball. Note that the Hölder elliptic 
regularity implies $u$ is $C^2$ (since the coefficients of the equation are Lipschitz) and this is what we need.

Next we observe that, as soon as $\abs{W_g}_g$ is positive, $\abs{W_g}_g^{\gamma}$ is smooth and
$$
L_g \abs{W_g}_g^{\gamma} = \Delta_g \abs{W_g}_g^{\gamma} - \gamma c(n) \abs{W}_g \abs{W_g}_g^{\gamma} \leq 0,
$$
so, with (\ref{invconf}),
$$
L_{\tilde{g}} ((1+u)^{-1} \abs{W_g}_g^{\gamma}) \leq 0,   
$$
which means 
$$
\Delta_{\tilde{g}} ((1+u)^{-1} \abs{W_g}_g^{\gamma}) \leq 0,
$$ 
outside a compact set.

Now, since $(M,\tilde{g})$ is quasi-isometric to $(M,g)$, it satisfies the doubling volume property as well as the scaled Poincaré inequality.
These properties are equivalent to the following two-sided gaussian estimate on the heat kernel $p_.(.,.)$ : for every $x,y$ in $M$, for 
every $t>0$, 
$$
\frac{c}{V(x,\sqrt{t})} \exp\left( -\frac{C d(x,y)^2}{t} \right) \leq p_t(x,y) \leq \frac{C}{V(x,\sqrt{t})} \exp\left( -\frac{c d(x,y)^2}{t} 
\right)
$$
(see \cite{SC2}, \cite{Grig}). As for large $R$, $V_{\tilde{g}}(o,R) \geq \tilde{A_o} R^\nu$, $\nu>2$, this in turn implies the existence of a positive Green function 
$G(.,.)$, which is simply $\int_0^\infty p_t(.,.) dt$) \cite{LY}. Using this formula and the upper bound on the heat kernel, we see that :
$$
G(o,x) = O(r_o(x)^{2-\nu})
$$
when $r_o(x)$ goes to infinity.
The maximum principle implies that for every point $x \in M \backslash B_g(o,R)$, 
$$
(1+u)^{-1} \abs{W_g}_g^{\gamma}(x) \leq \frac{\max_{S(o,R)} (1+u)^{-1} \abs{W_g}_g^{\gamma}}{\min_{S(o,R)} G(o,.)} G(o,x).
$$
We deduce 
$$
W = O(r^{\frac{2-\nu}{\gamma}}).
$$
\end{demo}

\begin{cor}
We consider a connected complete Ricci-flat manifold $M^n$, with $n \geq 4$. Assume there exists $o \in M$, 
$\nu > 4\frac{n-2}{n-1}$ and  $C_o >0$ such that 
$$
\forall R_2 \geq R_1 > 0, \  \frac{V(o,R_2)}{V(o,R_1)} \geq C_o \left( \frac{R_2}{R_1} \right)^\nu
$$
and the curvature satisfies 
$$
\int_M \abs{W}^{\frac{n}{2}} \rho_o(r_o) dvol < + \infty.
$$
Then 
$$
\sup_{S(o,R)} \abs{W} = O(r^{- \frac{(\nu-2)(n-1)}{n-3}}).
$$ 
\end{cor}

\begin{rem}
When $\nu=n=4$, we obtain the same decay as \cite{BKN}. 
\end{rem}

\begin{ex}
The Taub-NUT metric is a riemannian metric on $\Real^4$ introduced by Stephen Hawking in \cite{Haw} (see \cite{Leb} for a 
mathematical point of view). This is a Hyperkähler hence Ricci-flat metric with curvature decaying like $r^{-3}$ and volume growth 
like $r^3$. In this example, our theorem predicts the exact decay of the curvature.  
\end{ex}

\begin{ex}
Let us give another example, inspired from the famous Schwarzschild metric. We consider $\Real^2 \times \Sphere^{n-2}$, $n\geq 4$, endowed 
with the metric
$$
g = dr^2 + F(r)^2 dt^2 + G(r)^2 d\sigma^2.
$$   
$r$, $t$ are polar coordinates on the $\Real^2$ factor, $d\sigma^2$ is the standard metric on $\Sphere^{n-2}$, $F$ and $G$ are smooth
functions. Using the symmetries of this metric (see \cite{Bes}, \cite{Pet}), it is easy to obtain formulas for the curvature. And one sees 
that $g$ has vanishing Ricci tensor if and only if for some positive parameter $\gamma$, $G$ satisfies
$$
\left\lbrace
\begin{array}{l}
G'(r) = \sqrt{1 - \left(\frac{\gamma}{G}\right)^{n-3}}\\
G(0)=\gamma \\
G'(0)=0
\end{array}
\right.
$$
and 
$$
F(r) = \frac{2\gamma}{n-3} \sqrt{1 - \left(\frac{\gamma}{G}\right)^{n-3}}.
$$
$G$ increases from $\gamma$ to $\infty$ and $G \sim r$ at infinity ; $F$ increases from $0$ to $\frac{2\gamma}{n-3}$ and $F \sim r$ near
$0$. In particular, $g$ is $C^0$-close to the flat metric on $\Real^{n-1} \times \Sphere^1$ at infinity (the radius
of the circles at infinity are proportionnal to $\gamma$) and the distance to a fixed point in
this manifold behaves like the coordinate $r$ at infinity. Eventually, this provides on $\Real^2 \times \Sphere^{n-2}$, $n\geq 4$, a 
complete riemannian metric which is Ricci flat, has volume growing like $r^{n-1}$ and curvature decreasing like $r^{-(n-1)}$. This is what
our theorem predicted.   
\end{ex}


\vspace{2\baselineskip}



\end{document}